\documentclass[preprint,11pt]{elsarticle}

\usepackage{amssymb}
\usepackage{color,xcolor}
\usepackage{hyperref}
\usepackage{amsthm,epstopdf}
\usepackage{amsmath,amsfonts}
\usepackage{bookmark}
\usepackage[left=3.54cm,right=3.54cm,top=4.10cm,bottom=4.10cm]{geometry}
\biboptions{numbers,sort&compress}
\usepackage{lineno}

\modulolinenumbers[5]

\newtheorem{theorem}{Theorem}
\newtheorem{example}{Example}

\newtheorem{lemma}{Lemma}
\newtheorem{remark}{Remark}

\graphicspath{{figure_eps/}}
\usepackage{extarrows}
\usepackage{multicol}
\usepackage{array}
\usepackage{comment}
\newcommand{\ba}{\begin{array}}\newcommand{\ea}{\end{array}}
\newcommand{\be}{\begin{eqnarray}}\newcommand{\ee}{\end{eqnarray}}
\newcommand{\beq}{\begin{equation}}\newcommand{\eeq}{\end{equation}}
\newcommand{\bex}{\begin{eqnarray*}}
\newcommand{\eex}{\end{eqnarray*}}
\journal{}

\begin{document}
\title
{A novel spectral method for the subdiffusion equation}
\author[label1]{Chuanju Xu\corref{mycorrespondingauthor}}
\cortext[mycorrespondingauthor]{Corresponding author}
\ead{cjxu@xmu.edu.cn}
\author[label1]{Wei Zeng}

\address[label1]{School of Mathematical Sciences and Fujian Provincial Key Laboratory of Mathematical Modeling and High Performance Scientific Computing, Xiamen University, 361005 Xiamen, China}
%
%

\begin{frontmatter}
\begin{abstract}
In this paper, we design and analyze a novel
spectral method for the subdiffusion equation. As it has been known,
the solutions of this equation are usually singular near the initial time. Consequently,
direct application of the traditional high-order numerical methods is inefficient. We try to overcome this difficulty in a novel approach by combining variable transformation techniques with spectral methods.
The idea is to first use suitable variable transformation
to re-scale the underlying equation, then construct spectral methods for the re-scaled
equation. We establish a new variational framework based on the $\psi$-fractional Sobolev spaces. This allows us to prove the well-posedness of the associated variational problem. The proposed spectral method is based on the variational problem and
generalized Jacobi polynomials to approximate the re-scaled fractional
differential equation.
Our theoretical and numerical investigation show that the proposed method
is exponentially convergent for general right hand side functions, even though the exact solution has very limited regularity.
Implementation details are also provided, along with
a series of numerical examples to show the efficiency of the proposed method.
\end{abstract}

\end{frontmatter}


\section{Introduction}
\label{Sect1}

Fractional partial differential equations (FPDEs) appear in the
investigation of transport dynamics in complex systems which are
characterized  by the anomalous diffusion and nonexponential relaxation
patterns. Related equations of importance are the
space/time fractional diffusion equations, the fractional
advection-diffusion equation for anomalous diffusion
with sources and sinks, and the fractional Fokker-Planck equation
\cite{MK2,HW,HW2,BMK} for anomalous diffusion in an external field, etc. In
fact, it has been found that anomalous diffusion is
ubiquitous in physical and biological systems where trapping and
binding of particles can occur
\cite{BWZW,FBSBW,GW,SLSJ,SMWFC,SL,SM,MKW,AMYPL,Nig,Mai2}.
Anomalous diffusion deviates from
the standard Fichean description of Brownian motion, the main
character of which is that its mean squared displacement is a
nonlinear growth with respect to time, such as $\langle
x^2(t)\rangle\thicksim t^\alpha$.

The universality of anomalous diffusion phenomenon in physical and
biological experiments has led to an intensive investigation on the
fractional differential equations in recent years.
The time fractional diffusion equation (TFDE) considered in this
paper is of interest not only in its own right, but also in that it
constitutes the basic part in solving many other FPDEs.
The TFDE and related equations have been investigated in analytical
and numerical frames by a large number of authors,
see, e.g., \cite{SW06,LX07,GSZ12,GSZ14,LX15,Ali15,LX16,FY17,ZZK17,BH17,BH17SIAM,JZ17,ZZJZ17,YSZ17,ZTB17}.
Spectral methods have also been applied in solving the TFDE.
As is well known, any discretization including low-order approaches of a fractional
derivative has to take into account its non-local structure, which results in a full linear
system and a high storage requirement. Therefore it is very natural to consider a global
method, such as the spectral method, since the high accuracy of spectral methods may
significantly reduce the storage requirements. The first attempt in this direction was made by
Li and Xu in \cite{LIXU,LX10}. It was proven that the exponential convergence of
the proposed method is attainable for smooth solutions.

A main difficulty in numerically solving the TFDE comes from the fact that
the solutions of the TFDE are usually of low regularity, which lowers the accuracy of the above mentioned methods.
Some efforts have been made in developing and analyzing numerical methods for solutions of low regularity.
Among them, the modified time-stepping schemes are prominent, which can be roughly divided into two categories, i.e., piecewise polynomial interpolation based on a class of nonuniform grids and convolution quadrature with initial correction. Stynes et al. \cite{Stynes2017Error} proposed to use graded meshes in L1 scheme
for a reaction-diffusion problem, and an error analysis was given taking the starting time singularity into consideration. Later on, Liao et al. \cite{Liao2018Sharp, Liao2019A} gave a more general error analysis of L1 formula
on nonuniform grids based on a discrete fractional Gr\"onwall inequality.
Some researchers \cite{Lubich1996Nonsmooth, Cuesta2006Convolution, Jin2016Two, Jin2017Correction}
achieved optimal convergence rate by correcting the first several time steps.
Several other works focused
on spectral methods for non-smooth solutions of some related fractional equations, using
polyfractonomials \cite{ZK13,Zay15}, generalized Jacobi functions \cite{Chen2016Generalized},
M\"untz Jacobi polynomials \cite{Hou2017A, Hou2018Muntz},
and log orthogonal functions/generalized log orthogonal functions \cite{Chen2020Log,Chen2020A}.
Numerical experiments or theoretical analysis presented therein have shown
exponential convergence for non-smooth solutions having
specific singularity.

Unlike these existing numerical approaches,
in this paper we propose to first re-scale the time-fractional problems, then use the traditional approximations
to the re-scaled problems. Li et al. \cite{Li2021A} has tried this idea using a specific scaling function and
proposed two finite difference schemes based on the linear interpolation and quadratic interpolation.
The advantage of this approach is that the regularity of the re-scaled fractional operator can be much
higher than that of the original operator, which is more conducive to construct high-order schemes.
Below, we describe the main contributions of the paper and how the paper is organized.

Our first contribution is the development of the $\psi$-fractional Sobolev spaces presented
in Sect. \ref{Sect2}, which lays the foundation for the establishment of a new variational framework
in Sect. \ref{Sect3}. In detail, we introduce the concept of $\psi$--fractional operators,
and propose $\psi$--fractional Sobolev spaces on this basis and prove the equivalence of related norms of $\psi$--fractional Sobolev spaces.

The second contribution is to propose a new Galerkin spectral method based on the generalized Jacobi polynomials
(GJPs) under the new variational framework introduced in Sect. \ref{Sect3}.
The well-posedness of the weak problem is proved in the $\psi$-fractional Sobolev spaces, together with
error estimation established in the non-uniform Jacobi-weighted norm.
Moreover, it's shown that the new approach is as efficient as the M\"untz spectral method \cite{Hou2017A, Hou2018Muntz} by using suitable scaling functions.
The novel approach not only provides a theoretical support of the M\"untz spectral method
but also gives a guideline for the selection of the scaling parameters.

Finally, the proposed approach is applied to the time fractional subdiffusion equations in Sect. \ref{Sect4}.
A space-time Galerkin spectral method is developed based on the re-scaled weak formulation
and a combination of temporal GJPs and spatial Legendre polynomials.
In Sect. \ref{Sect5}, we present some numerical tests to confirm the theoretical findings.
Some concluding remarks are given in Sect. \ref{Sect6}.

\section{Functional Spaces}\label{Sect2}

In order to develop the re-scaling method for fractional differential equations,
we need some preparations, mainly including an introduction of the $\psi$-fractional Sobolev spaces
and establishment of the associated variational framework.
Throughout this paper, let $c$ stand for a generic positive constant
independent of any functions and of any discretization parameters.
In what follows, we use the expression $A\lesssim B$
(respectively, $A\gtrsim B$) to mean that $A\leqslant cB$ (respectively, $A\geqslant cB$).

The first part of this section is devoted to introducing the $\psi$--fractional integrals,
derivatives and a crucial variable transformation.

\subsection{$\psi$--fractional operators and variable transformation}\label{Sec2_subsec1}
We recall some definitions of $\psi$--fractional integrals and
$\psi$--fractional derivatives; see Kilbas et al. \cite[Sect. 2.5]{Kilbas2006Theory} or Samko et al.
\cite[Sect. 18.2]{Samko1993Fractional}.
Let $\Gamma(\cdot)$ denote the Gamma function.
For any positive integer $n$  and real number $\delta$, $n-1\leq \delta<n$,
$v$ is an integrable function in the bounded interval $[a_\psi, b_\psi]$ with respect to the function $\psi:[a_\psi,b_\psi]\to \mathbb{R}$ that is increasing and differentiable
such that $\psi^{\prime}(t)\neq 0$.
The $\psi$--fractional integral, $\psi$--Caputo derivative, and $\psi$--Riemann--Liouville derivative of order
$\delta $ of $v$ are respectively defined as follows:
 $\forall~t\in [a_\psi,b_\psi],$
\begin{equation}\label{lRI}
\mbox{left $\psi$--fractional integral:}\quad
I^{\!\delta,\psi}_t v(t) := \frac{1}{\Gamma(\delta)}\int_{a_\psi}^t\psi^{\prime}(z)(\psi(t)-\psi(z))^{\delta-1}v(z)dz,\tag{I1}
\end{equation}
\begin{equation}\label{rRI}
\mbox{right $\psi$--fractional integral:}\quad
_tI^{\!\delta,\psi} v(t) := \frac{1}{\Gamma(\delta)}\int_t^{b_\psi}\psi^{\prime}(z)(\psi(z)-\psi(t))^{\delta-1}v(z)dz,\tag{I2}
\end{equation}
\begin{equation}\label{lC}
\mbox{left $\psi$--Caputo derivative:}\quad
^C\!D^{\!\delta,\psi}_tv(t) := I^{n-\delta,\psi}_t \left(\frac{1}{\psi^{\prime}(t)}\frac{d}{dt}\right)^nv(t), \tag{D1}
\end{equation}
\begin{equation}\label{rCD}
\mbox{right $\psi$--Caputo derivative:}\quad
^C_t\!D^{\!\delta,\psi}v(t) :=\,_t\!I^{n-\delta,\psi} \left(-\frac{1}{\psi^{\prime}(t)}\frac{d}{dt}\right)^nv(t),\tag{D2}
\end{equation}
\begin{equation}\label{lRD}
\mbox{left $\psi$--Riemann--Liouville derivative:}\quad
\!D^{\!\delta,\psi}_tv(t) = \left(\frac{1}{\psi^{\prime}(t)}\frac{d}{dt}\right)^n\,I^{n-\delta,\psi}_t v(t),\tag{D3}
\end{equation}
\begin{equation}\label{rRD}
\mbox{right $\psi$--Riemann--Liouville derivative:}~
_t\!D^{\!\delta,\psi}v(t) := \left(-\frac{1}{\psi^{\prime}(t)}\frac{d}{dt}\right)^n \,_t\!I^{n-\delta,\psi} v(t)\tag{D4}.
\end{equation}
When $\psi(t)=t,$ the above definitions degenerate into the classical fractional integral, Caputo derivative and Riemann--Liouville derivative; see \cite{Oldham1974The, Podlubny1999Fractional}. In particular,
when $\delta\in (0,1),$ $\psi$--Caputo fractional derivative $^C\!D^{\!\delta,\psi}_tv(t)$ becomes
$^C\!D^{\delta}_tv(t),$ where
\begin{equation}\label{classic_Caputo}
	^C\!D^{\delta}_tv(t)
	:=\frac{1}{\Gamma(1-\delta)}\int_0^t(t-z)^{-\delta}v^{\prime}(z)dz.
\end{equation}
On the contrary, by a change of variable $\psi(t)$, the classical Caputo fractional derivative can be turned into a class of $\psi$--Caputo fractional derivative.
For example, a direct calculation gives
\begin{align}\label{var_trans}
	\begin{split}
	^C\!D^\delta_s u(s)=&\frac{1}{\Gamma(1-\delta)}\int_0^s(s -x)^{-\delta}u^{\prime}(x)dx,\quad s\in (0,T]\\
	\xlongequal[]{s=t^{1/\gamma}}&\frac{1}{\Gamma(1-\delta)}\int_0^{t^{1/\gamma}}\big(t^{1/\gamma} -x)^{-\delta}u^{\prime}(x)dx,\quad t\in (0,T^{\gamma}]\\
	\xlongequal[]{x=z^{1/\gamma}}&\frac{1}{\Gamma(1-\delta)}\int_0^t\big(t^{1/\gamma} -z^{1/\gamma}\big)^{-\delta}du(z^{1/\gamma}),\quad t\in (0,T^{\gamma}].
	\end{split}
\end{align}
Let $\psi(t)=t^{1/\gamma}$, $v(t):=u(\psi(t))$ for $t\in (0,T^{\gamma}]$. Then the new fractional derivative
$^C_*\!D_t^{\delta}v(t)$, defined by
\begin{align}\label{new_Caputo}
	^C_*\!D_t^{\delta}v(t):=\frac{1}{\Gamma(1-\delta)}\int_0^t\big(t^{1/\gamma} -z^{1/\gamma}\big)^{-\delta}v^{\prime}(z)dz,\quad t\in (0,T^{\gamma}]
\end{align}
can be regarded as a class of $\psi$--Caputo fractional derivative of $v(t)$ with $\psi(t)=t^{1/\gamma}$.

The following results about $\psi$--fractional operators are frequently used; see \cite{Almeida2017A, Almeida2018Fractional}.

It is noted that the $\psi$--Riemann--Liouville fractional derivative and $\psi$--Caputo fractional derivative of $v$ have the following relationship
\be\label{C-RL}
^C\!D^{\delta,\psi}_tv(t)=D^{\delta,\psi}_t[v(t)-v(a_\psi)],\quad \delta\in(0,1).
\ee
And left $\psi$--Riemann--Liouville fractional derivative and integral of order $\delta$ satisfy
\be\label{deri_jifen}
D^{\!\delta,\psi}_t I^{\!\delta,\psi}_tv(t) = v(t).
\ee

With the above notations and properties, we are in a position to introduce the $\psi$-fractional Sobolev spaces.
\subsection{$\psi$-fractional Sobolev spaces}
We begin with some additional notations. Let $0<\gamma\leq 1, I=(a,b)$.
The function $\psi^{-1}(\cdot)$ denotes the inverse function of $\psi(\cdot)$.
Let $\Lambda=(a_\psi, b_\psi):=(\psi^{-1}(a),\psi^{-1}(b))$.
Thus if $s\in I$, then $t=\psi^{-1}(s)\in \Lambda$.

Define the space
$$L^2_{\psi}(\Lambda)=\left\{v:\Lambda\rightarrow \mathbb{R} \big| v~\mathrm{is~measurable~and}~ \int_{\Lambda}|v(t)|^2\psi^{\prime}(t)dt<\infty\right\}.$$
It can be easily seen that $L^2_{\psi}(\Lambda)$ is a Hilbert space with respect to the scalar product
\be\label{innerp}
( v,w)_{L^2_{\psi}(\Lambda)}=\int_{\Lambda}v(t)w(t)\psi^{\prime}(t)dt,~~
\forall v,w\in L^2_{\psi}(\Lambda).
\ee
The norm in $L^2_{\psi}(\Lambda)$ induced by the scalar product $(\cdot,\cdot)_{L^2_{\psi}(\Lambda)}$
is defined by
$$\|v\|_{L^2_{\psi}(\Lambda)}=\sqrt{( v,v)_{L^2_{\psi}(\Lambda)}}=\left(\int_{\Lambda}|v(t)|^2\psi^{\prime}(t)dt\right)^{1/2},~~\forall v\in L^2_{\psi}(\Lambda).$$
In particular, for $\psi(t)=t,$ the space $L^2_{\psi}(\Lambda)$ is reduced to the classical space $L^2(\Lambda).$ Let us denote by $(\cdot,\cdot)_{L^2(\Lambda)}$ and $\|\cdot,\cdot\|_{L^2(\Lambda)}$
the inner product and norm in $L^2(\Lambda),$ respectively.

We now introduce the $\psi$--fractional Sobolev spaces.
Let $\mathcal{F}(v)$ denote the Fourier transform of $v$, $\tilde{v}(\cdot):=v(\psi^{-1}(\cdot))$. Define the space
\be\label{Hpsi}
H^{\!\delta,\psi}(\mathbb{R}):=\left\{v\big|v\in L^2_{\psi}(\mathbb{R});\ (1+|\xi|^{\delta})\mathcal{F}(\tilde{v})(\xi)\in L^2(\mathbb{R})\right\},\ \delta\geq 0,
\ee
endowed with the semi-norm and norm
$$\left|v\right|_{H^{\!\delta,\psi}(\mathbb{R})}=\left\||\xi|^{\delta}\mathcal{F}(\tilde{v})(\xi)\right\|_{L^2(\mathbb{R})},$$
$$\left\|v\right\|_{H^{\!\delta,\psi}(\mathbb{R})}=\left(\left\|v\right\|^2_{L^2_{\psi}(\mathbb{R})}+\left|v\right|^2_{H^{\!\delta,\psi}(\mathbb{R})}\right)^{1/2},$$
respectively.
Note that $\mathcal{F}(\tilde{v})$ rather than $\mathcal{F}({v})$ was used in the definition \eqref{Hpsi}.

The $\psi$--fractional Sobolev space for the bounded domain $\Lambda$ is defined by
$$H^{\!\delta,\psi}(\Lambda):=\left\{v\in L^2_{\psi}(\Lambda)\big|\exists\ v_{e}\in H^{\!\delta,\psi}(\mathbb{R})\ \mathrm{such}\ \mathrm{that}\ v_{e}|_{\Lambda}=v\right\},$$
equipped with the norm
$$\|v\|_{H^{\!\delta,\psi}(\Lambda)}=\inf\limits_{v_{e}\in H^{\!\delta,\psi}(\mathbb{R}),\ v_{e}|_{\Lambda}=v}\left\|v_{e}\right\|_{H^{\!\delta,\psi}(\mathbb{R})}.$$
It is readily seen that $H^{\!\delta,\psi}(\Lambda)$ degenerates into the classic Sobolev space
$H^{\delta}(\Lambda)$ when $\psi(t)=t.$

We define
$$^{L\!}H^{\!\delta,\psi}(\Lambda):=\left\{v\big|\|v\|_{^{L\!}H^{\!\delta,\psi}(\Lambda)}<\infty\right\},$$
where $|\cdot|_{^{L\!}H^{\!\delta,\psi}(\Lambda)}$ is the norm:
\begin{equation*}
\left\|v\right\|_{^{L\!}H^{\!\delta,\psi}(\Lambda)}=\left(\|v\|^2_{L^2_{\psi}(\Lambda)}+|v|^2_{^{L\!}H^{\!\delta,\psi}(\Lambda)}\right)^{1/2},\ \ |v|_{^{L\!}H^{\!\delta,\psi}(\Lambda)}=\left\|D^{\!\delta,\psi}_tv\right\|_{L^2_{\psi}(\Lambda)}.
\end{equation*}
Similarly, we define
$$^{R\!}H^{\!\delta,\psi}(\Lambda):=\left\{v\big|\|v\|_{^{R\!}H^{\!\delta,\psi}(\Lambda)}<\infty\right\},$$
with
\begin{equation*}
\|v\|_{^{R\!}H^{\!\delta,\psi}(\Lambda)}=\left(\|v\|^2_{L^2_{\psi}(\Lambda)}+|v|^2_{^{R\!}H^{\!\delta,\psi}(\Lambda)}\right)^{1/2},
\ \ |v|_{^{R\!}H^{\!\delta,\psi}(\Lambda)}=\left\|{}_tD^{\!\delta,\psi}v\right\|_{L^2_{\psi}(\Lambda)};
\end{equation*}
and
$$^{c\!}H^{\!\delta,\psi}(\Lambda):=\left\{v\big|\|v\|_{^{c\!}H^{\!\delta,\psi}(\Lambda)}<\infty\right\},$$ with
\begin{equation}
	\|v\|_{^{c\!}H^{\!\delta,\psi}(\Lambda)}
	\!=\!\left(\|v\|^2_{L^2_{\psi}(\Lambda)}+|v|^2_{^{c\!}H^{\!\delta,\psi}(\Lambda)}\right)^{1/2}\!,
	\ |v|_{^{c\!}H^{\!\delta,\psi}(\Lambda)}
	\!=\!\left|(D^{\!\delta,\psi}_tv,{}_tD^{\!\delta,\psi}v)_{L^2_{\psi}(\Lambda)}\right|^{1/2}\!.
\end{equation}

Let $C^{\infty}_0(\Lambda)$ is the space of smooth functions with compact support in $\Lambda.$
Let $^{L\!}H^{\!\delta,\psi}_0(\Lambda)$, $^{R\!}H^{\!\delta,\psi}_0(\Lambda)$, $^{c\!}H_0^{\!\delta,\psi}(\Lambda)$,
and $H^{\!\delta,\psi}_0(\Lambda)$ be the closures of $C^{\infty}_0(\Lambda)$ with respect to the norms $\|v\|_{^{L\!}H^{\!\delta,\psi}(\Lambda)}$, $\|v\|_{^{R\!}H^{\!\delta,\psi}(\Lambda)}$, $\|v\|_{^{c\!}H^{\!\delta,\psi}(\Lambda)}$ and $\|v\|_{H^{\!\delta,\psi}(\Lambda)},$ respectively.
Besides, let ${}_0H^{\!\delta,\psi}(\Lambda)$ denote the closure of ${}_0C^{\infty}(\Lambda)$ with respect to $\|\cdot\|_{H^{\!\delta,\psi}(\Lambda)},$ where ${}_0C^{\infty}(\Lambda)$ is the space of smooth functions with compact support in $(a_\psi,b_\psi].$

Next we give some crucial lemmas, especially the equivalence results of the related norms
of the $\psi$-fractional Sobolev spaces. These results play a key role
in the subsequent analysis, including the well-posedness analysis and error estimation of
the numerical methods to be constructed.

\subsection{Some useful Lemmas}
Define the convolution of the functions $h_1(t)$ and $h_2(t)$ as follows:
$$h_1(t)\ast h_2(t):=\int^{+\infty}_{-\infty}h_1(t-\tau)h_2(\tau)d\tau
=\int^{+\infty}_{-\infty}h_1(\tau)h_2(t-\tau)d\tau,$$
where $h_1(t),$ $h_2(t)\in (-\infty,+\infty).$

It is known \cite{Podlubny1999Fractional} that
if the Fourier transform of $h_1(t)$ and $h_2(t)$ exists, then
\begin{equation}\label{Four_convo}
\mathcal{F}\{h_1(t)\ast h_2(t);~\xi\}=H_1(\xi)H_2(\xi),
\end{equation}
where $H_1(\xi)=\mathcal{F}\{h_1(t);~\xi\},$ $H_2(\xi)=\mathcal{F}\{h_2(t);~\xi\}.$

Then we can define Fourier transform of $\psi$- fractional derivatives on the above basis.
\begin{lemma}\label{Four_Frac}
	~(Fourier transform of $\psi$--fractional derivatives)~Let $v\in C^{\infty}_{0}(\mathbb{R}),$ $0<\delta<1.$
	Assume $\psi(\infty)=\infty$. Then
	\begin{equation}\label{Fourier_Frac}
		\mathcal{F}(_{-\infty}D^{\!\delta,\psi}_tv(t))=(i\xi)^{\delta}\mathcal{F}(\tilde{v})(\xi),~~~
		\mathcal{F}(_{t}D^{\!\delta,\psi}_{+\infty}v(t))=(-i\xi)^{\delta}\mathcal{F}(\tilde{v})(\xi).
	\end{equation}
\end{lemma}

\begin{proof}
	We first evaluate the Fourier transform of the $\psi$--fractional integral $_{-\infty}I^{\!\delta,\psi}_tv.$
The Laplace transform $\mathcal{L}$ of $s^{\delta-1}$ reads
\begin{equation}\label{integral_L}
\mathcal{L}\{s^{\delta-1};\tau\}=\int^{\infty}_0 s^{\delta-1}e^{-\tau s}ds=\Gamma(\delta)\tau^{-\delta}.
\end{equation}
Note that the above integral makes sense for all $\delta>0$ by the Dirichlet theorem.
Let $h(s)$ be the function
\begin{align*}\label{hf}
	\begin{split}
		h(s)=
		\begin{cases}
			\frac{s^{\delta-1}}{\Gamma(\delta)},\ &s>0,\\
			0,                                          \ &s\leq 0.
		\end{cases}
	\end{split}
\end{align*}
Then a direct calculation using \eqref{integral_L} shows
\begin{equation}\label{FTh}
\mathcal{F}\{h(s);\xi\}
=\int^{\infty}_{-\infty}h(s)e^{-i\xi s}ds
=\frac{1}{\Gamma(\delta)}\int^{\infty}_0 s^{\delta-1}e^{-i\xi s}ds
=\frac{\mathcal{L}\{s^{\delta-1};i\xi\}}{\Gamma(\delta)}
=(i\xi)^{-\delta}.
\end{equation}

Let $s=\psi(t),$ then the $\psi$--fractional integral $_{-\infty}I^{\!\delta,\psi}_tv$
can be expressed as a convolution of the functions $h(\psi(t))$ and $\tilde{v}(\psi(t)),$ i.e.,
	\begin{align*}
		\begin{split}
		_{-\infty}I^{\!\delta,\psi}_tv(t)&=\frac{1}{\Gamma(\delta)}\int_{-\infty}^{\psi(t)}(\psi(t)-\psi(\tau))^{\delta-1}\tilde{v}(\psi(\tau))d\psi(\tau)\\
		&=h(s)\ast \tilde{v}(s),
			\end{split}
\end{align*}
where $\tilde{v}(\cdot):=v(\psi^{-1}(\cdot)).$
Thus, it follows from \eqref{Four_convo} and \eqref{FTh}:
	\begin{equation*}
		\mathcal{F}\{_{-\infty}I^{\!\delta,\psi}_tv(t);\xi\}=(i\xi)^{-\delta}\mathcal{F}(\tilde{v})(\xi).
	\end{equation*}
Next, we calculate the Fourier transform of the $\psi$--fractional derivatives.
Note that
\bex
_{-\infty}D^{\!\delta,\psi}_tv(t)
=_{-\infty}\!I^{n-\delta,\psi}_t\left(\frac{1}{\psi^{\prime}(t)}\frac{d}{dt}\right)^nv(t),
\eex	
we have
\begin{align*}
		\begin{split}
			\mathcal{F}\{_{-\infty}D^{\!\delta,\psi}_t v(t)\}
			&=(i\xi)^{\delta-n}\mathcal{F}\left\{\left(\frac{d}{d\psi}\right)^n\tilde{v}(\psi(t));\xi\right\}\\
			&=(i\xi)^{\delta-n}(i\xi)^{n}\mathcal{F}(\tilde{v})(\xi)\\
			&=(i\xi)^{\delta}\mathcal{F}(\tilde{v})(\xi).
		\end{split}
	\end{align*}
The second equality in \eqref{Fourier_Frac} can be proved in a similar way.
\end{proof}

With the help of the Fourier transform of the $\psi$--fractional derivatives,
we can derive the following equivalence result for the $\psi$--fractional Sobolev spaces
on the whole line $\mathbb{R}$.
\begin{lemma}\label{equal_R}
	Let $\delta>0,~\delta\neq n-1/2,~n\in \mathbb{N}.$ Then the spaces $^{L\!}H^{\!\delta,\psi}(\mathbb{R}),$ $^{R\!}H^{\!\delta,\psi}(\mathbb{R}),$  $^{c\!}H^{\!\delta,\psi}(\mathbb{R})$, and $H^{\!\delta,\psi}(\mathbb{R})$
are equal to each other with equivalent semi-norms and norms.
\end{lemma}
\begin{proof}
	The proof will be divided into three steps.
	
	Step 1: the equivalence of the spaces $^{L\!}H^{\!\delta,\psi}(\mathbb{R})$ and $H^{\!\delta,\psi}(\mathbb{R}).$
	
	For a function $v\in {}^{L\!}H^{\!\delta,\psi}(\mathbb{R}),$
	we have $D^{\!\delta,\psi}_tv\in L^2_{\psi}(\mathbb{R}).$
	Using Lemma \ref{Four_Frac} and Plancherel's theorem gives
	\begin{equation*}
		\int_{\mathbb{R}}\left|\xi\right|^{2\delta}\left|\mathcal{F}(\tilde{v})(\xi)\right|^{2} d\xi=\int_{\mathbb{R}}\left|D^{\!\delta,\psi}_tv\right|^{2}\psi^{\prime}(t)dt.
	\end{equation*}
	Thus,
	\begin{equation*}
		\left\|\left|\xi\right|^{\delta}\mathcal{F}(\tilde{v})(\xi)\right\|_{L^{2}(\mathbb{R})}
		=\left|v\right|_{^{L\!}H^{\!\delta,\psi}(\mathbb{R})}.
	\end{equation*}
The desired result follows immediately from the above equality and the definition of $H^{\!\delta,\psi}(\mathbb{R})$.
	
	Step 2: the equivalence of the spaces $^{L\!}H^{\!\delta,\psi}(\mathbb{R})$ and $^{R\!}H^{\!\delta,\psi}(\mathbb{R}).$
	
	Again, using the results of Lemma \ref{Four_Frac} and Plancherel's theorem, we have
	\begin{equation}\label{Left_semi}
		\left|v\right|^2_{^{L\!}H^{\!\delta,\psi}(\mathbb{R})}
		=\int_{\mathbb{R}}\left|(i\xi)^{\delta}\mathcal{F}(\tilde{v})(\xi)\right|^{2} d\xi.
	\end{equation}
	Similarly,
	\begin{equation}\label{Right_semi}
		\left|v\right|^2_{^{R\!}H^{\!\delta,\psi}(\mathbb{R})}
		=\int_{\mathbb{R}}\left|(-i\xi)^{\delta}\mathcal{F}(\tilde{v})(\xi)\right|^{2} d\xi.
	\end{equation}
	Note that $|(i\xi)^{\delta}|=|(-i\xi)^{\delta}|.$ Thus the semi-norms
	$|v|_{^{L\!}H^{\!\delta,\psi}(\mathbb{R})}$
	and $|v|_{^{R\!}H^{\!\delta,\psi}(\mathbb{R})}$, consequently the norms
	$\|v\|_{^{L\!}H^{\!\delta,\psi}(\mathbb{R})}$
	and $\|v\|_{^{R\!}H^{\!\delta,\psi}(\mathbb{R})}$, are equivalent.
	
	Step 3: the equivalence of the spaces $^{c\!}H^{\!\delta,\psi}(\mathbb{R})$ and $^{L\!}H^{\!\delta,\psi}(\mathbb{R}).$
	
	Analogous to \cite[Lemma 2.4]{Ervin2006Variational}, with the help of some related properties of the
	Fourier transform, we obtain
	\begin{equation}\label{c_cos}
		\left(D^{\!\delta,\psi}_tv,~_tD^{\!\delta,\psi}v\right)_{L^2_{\psi}(\mathbb{R})}
		=\mathrm{cos}(\pi\delta)\left\|D^{\!\delta,\psi}_tv\right\|^2_{L^2_{\psi}(\mathbb{R})}.
	\end{equation}
That is,
	\begin{equation*}
		\left|v\right|^2_{^{c\!}H^{\!\delta,\psi}(\mathbb{R})}
		=|\mathrm{cos}(\pi\delta)|\left|v\right|^2_{^{L\!}H^{\!\delta,\psi}(\mathbb{R})}.
	\end{equation*}
Thus the semi-norms of $^{c\!}H^{\!\delta,\psi}(\mathbb{R})$ and $^{L\!}H^{\!\delta,\psi}(\mathbb{R})$ are
equivalent. So are their norms, which implies the equivalence of
$^{c\!}H^{\!\delta,\psi}(\mathbb{R})$ and $^{L\!}H^{\!\delta,\psi}(\mathbb{R}).$

We conclude by combining Step 1--Step 3.
\end{proof}

The equivalence of different $\psi$--fractional spaces on the bounded interval $\Lambda$ are established below.
\begin{lemma}\label{lem_equi}
	Let $\delta>0,~\delta\neq n-1/2,~n\in \mathbb{N}.$ Then the spaces $^{L\!}H^{\!\delta,\psi}_0(\Lambda),$ $^{R\!}H^{\!\delta,\psi}_0(\Lambda),$  $^{c\!}H^{\!\delta,\psi}_0(\Lambda)$, and $H^{\!\delta,\psi}_0(\Lambda)$
are equal to each other with equivalent semi-norms and norms.
\end{lemma}
\begin{proof}
The proof is splitted into two steps.
	
	Step 1: the equivalence of the spaces $^{c\!}H^{\!\delta,\psi}_0(\Lambda)$ and $H^{\!\delta,\psi}_0(\Lambda).$
	
	For $v\in C^{\infty}_0(\Lambda)$, let $v_e$ be the extension of $v$ by zero outside of $\Lambda.$ Then
	\begin{align*}
		\begin{split}
			&\mathrm{supp}\left(v_e\right)\subset \Lambda,\\
			&\mathrm{supp}\left(D^{\!\delta,\psi}_t v_e\right)\subset (a_\psi,\infty),\\
			&\mathrm{supp}\left(_tD^{\!\delta,\psi} v_e\right)\subset (-\infty,b_\psi).
		\end{split}
	\end{align*}
	Thus,
	\begin{equation*}
		\mathrm{supp}\left(D^{\!\delta,\psi}_t v_e\ {_tD^{\!\delta,\psi}}v_e\right)\subset \Lambda,
	\end{equation*}
from which it follows
	\begin{equation*}
		\left|v\right|_{^{c\!}H^{\!\delta,\psi}_0(\Lambda)}
		=\left|v_e\right|_{^{c\!}H^{\!\delta,\psi}_0(\mathbb{R})}.
	\end{equation*}
On the other side, we have
	\begin{equation*}
		\left|v\right|_{H^{\!\delta,\psi}_0(\Lambda)}
		=\left|v_e\right|_{H^{\!\delta,\psi}_0(\mathbb{R})}.
	\end{equation*}
Then the semi-norm equivalence of ${^{c\!}H^{\!\delta,\psi}_0(\mathbb{R})}$ and
${H^{\!\delta,\psi}_0(\mathbb{R})}$, proved in Lemma \ref{equal_R}, gives
	\begin{equation*}
		\left|v\right|_{^{c\!}H^{\!\delta,\psi}_0(\Lambda)}
		=\left|v\right|_{H^{\!\delta,\psi}_0(\Lambda)}.
	\end{equation*}
Thus the spaces $^{c\!}H^{\!\delta,\psi}_0(\Lambda)$ and
	$H^{\!\delta,\psi}_0(\Lambda)$ are equal with equivalent norms.
	
	Step 2: the equivalence of the spaces $^{L\!}H^{\!\delta,\psi}_0(\Lambda),$ $^{R\!}H^{\!\delta,\psi}_0(\Lambda)$, and $H^{\!\delta,\psi}_0(\Lambda).$
	
It follows from \eqref{Left_semi} and the definition of $|\cdot|_{H^{\!\delta,\psi}_0(\Lambda)}$:
	\begin{equation*}
		\left\|D^{\!\delta,\psi}_t v\right\|_{L^2_{\psi}(\Lambda)}=\left|v\right|_{^{L\!}H^{\!\delta,\psi}(\Lambda)}
		\leq \left|v_e\right|_{^{L\!}H^{\!\delta,\psi}(\mathbb{R})}
		=\left|v_e\right|_{H^{\!\delta,\psi}(\mathbb{R})}
		=\left|v\right|_{H^{\!\delta,\psi}_0(\Lambda)}.
	\end{equation*}
This gives
	\begin{equation*}
		H^{\!\delta,\psi}_0(\Lambda)\subseteq\ ^{L\!}H^{\!\delta,\psi}_0(\Lambda).
	\end{equation*}
	Combining the result proved in Step 1 and Young's inequality, we obtain
	\begin{align*}
		\begin{split}
			\left|v\right|^2_{H^{\!\delta,\psi}_0(\Lambda)}
			&\leq c\left|v\right|^2_{^{c\!}H^{\!\delta,\psi}_0(\Lambda)}\\
			&=c\left|\int_{\Lambda}D^{\!\delta,\psi}_tv\ _tD^{\!\delta,\psi}v\psi^{\prime}(t)dt\right|\\
			&\leq \frac{c}{4\epsilon}\left\|D^{\!\delta,\psi}_tv\right\|^2_{L^2_{\psi}(\Lambda)}
			+\epsilon c\left\|_tD^{\!\delta,\psi}v\right\|^2_{L^2_{\psi}(\Lambda)}.
		\end{split}
	\end{align*}
Furthermore, it follows from \eqref{Right_semi} and the definition of
$|\cdot|_{H^{\!\delta,\psi}_0(\Lambda)}$:
	\begin{equation*}
		\left\|_tD^{\!\delta,\psi} v\right\|_{L^2_{\psi}(\Lambda)}=\left|v\right|_{^{R\!}H^{\!\delta,\psi}(\Lambda)}
		\leq \left|v_e\right|_{^{R\!}H^{\!\delta,\psi}(\mathbb{R})}
		=\left|v_e\right|_{H^{\!\delta,\psi}(\mathbb{R})}
		=\left|v\right|_{H^{\!\delta,\psi}_0(\Lambda)}.
	\end{equation*}
Combining the last two inequalities gives
	\begin{equation*}
		\left|v\right|^2_{H^{\!\delta,\psi}_0(\Lambda)}\leq
		\frac{c}{2\epsilon}\left|v\right|^2_{^{L\!}H^{\!\delta,\psi}_0(\Lambda)}
		+\frac{\epsilon c}{2}\left|v\right|^2_{H^{\!\delta,\psi}_0(\Lambda)}.
	\end{equation*}
	Taking $\epsilon=1/c$ in the above inequality yields
	\begin{equation*}
		\left|v\right|^2_{H^{\!\delta,\psi}_0(\Lambda)}
		\leq c^2\left|v\right|^2_{^{L\!}H^{\!\delta,\psi}_0(\Lambda)}.
	\end{equation*}
This gives
	\begin{equation*}
		^{L\!}H^{\!\delta,\psi}_0(\Lambda)\subseteq\ H^{\!\delta,\psi}_0(\Lambda).
	\end{equation*}
This ends the proof of the semi-norm equivalence of the spaces
$^{L\!}H^{\!\delta,\psi}_0(\Lambda)$ and $H^{\!\delta,\psi}_0(\Lambda)$, and thus the equivalence of the spaces themselves.
In a similar way, we can prove	
the equivalence of the spaces $^{R\!}H^{\!\delta,\psi}_0(\Lambda)$ and
$H^{\!\delta,\psi}_0(\Lambda)$. The proof is completed.
\end{proof}

Now we turn to derive some Poincar\'e-Friedrichs-type inequalities for the functions
in $\psi$--fractional spaces.
The following mapping properties are useful.

\begin{lemma}\label{map_lem}(Mapping properties) All the following mappings are
bounded linear operator.
\begin{enumerate}
		\renewcommand\theenumi{\roman{enumi}}
		\renewcommand{\labelenumi}{$\mathrm{\left(\theenumi\right)}$}
	\item $I^{\!\delta,\psi}_t:~L^2_{\psi}(\Lambda)\rightarrow L^2_{\psi}(\Lambda)$.
	\item $I^{\!\delta,\psi}_t:~L^2_{\psi}(\Lambda)\rightarrow~^{L\!}H^{\!\delta,\psi}(\Lambda)$.
	\item $D^{\!\delta,\psi}_t:~^{L\!}H^{\!\delta,\psi}(\Lambda)\rightarrow L^2_{\psi}(\Lambda)$.
	\item $_{t}I^{\!\delta,\psi}:~L^2_{\psi}(\Lambda)\rightarrow L^2_{\psi}(\Lambda)$.
	\item $_{t}I^{\!\delta,\psi}:~L^2_{\psi}(\Lambda)\rightarrow\ ^{R\!}H^{\!\delta,\psi}(\Lambda)$.
	\item $_{t}D^{\!\delta,\psi}:~^{R\!}H^{\!\delta,\psi}(\Lambda)\rightarrow L^2_{\psi}(\Lambda)$.
	\end{enumerate}
\end{lemma}
\begin{proof}
	\begin{enumerate}
		\renewcommand\theenumi{\roman{enumi}}
		\renewcommand{\labelenumi}{$\mathrm{\left(\theenumi\right)}$}
		\item By using \cite[Lemma 2.2]{Almeida2019A}, we have
		\begin{equation*}
			\left\|I^{\!\delta,\psi}_t v\right\|_{L^2_{\psi}(\Lambda)}\leq \frac{1}{\Gamma(\delta)\sqrt{2\delta-1}}\left\|v\right\|_{L^2_{\psi}(\Lambda)},
		\end{equation*}
which means $I^{\!\delta,\psi}_t$ is a bounded linear operator from
$L^2_{\psi}(\Lambda)$ to $L^2_{\psi}(\Lambda)$.
		\item Combining \eqref{deri_jifen} and the definition of $^{L\!}H^{\!\delta,\psi}(\Lambda),$ then using (i), one obtains
		\begin{align*}
			\begin{split}
				\left\|I^{\!\delta,\psi}_t v\right\|_{^{L\!}H^{\!\delta,\psi}(\Lambda)}
				&=\left(\left\|I^{\!\delta,\psi}_tv\right\|^2_{L^2_{\psi}(\Lambda)}
				+\left\|D^{\!\delta,\psi}_tI^{\!\delta,\psi}_tv\right\|^2_{L^2_{\psi}(\Lambda)}\right)^{1/2}\\
				&\lesssim \left\|v\right\|_{L^2_{\psi}(\Lambda)}.
			\end{split}
		\end{align*}
This proves that $I^{\!\delta,\psi}_t$ is a bounded linear operator from
$L^2_{\psi}(\Lambda)$ to $^{L\!}H^{\!\delta,\psi}(\Lambda)$.
		\item It follows from the definition of the norm $\|\cdot\|_{^{L\!}H^{\!\delta,\psi}(\Lambda)}$:
		\begin{align*}
			\begin{split}
				\left\|D^{\!\delta,\psi}_t v\right\|_{L^2_{\psi}(\Lambda)}
				&\leq \left(\left\|v\right\|^2_{L^2_{\psi}(\Lambda)}
				+\left\|D^{\!\delta,\psi}_tv\right\|^2_{L^2_{\psi}(\Lambda)}\right)^{1/2}\\
				&=: \left\|v\right\|_{^{L\!}H^{\!\delta,\psi}(\Lambda)}.
			\end{split}
		\end{align*}
This shows that $D^{\!\delta,\psi}_t$ is a bounded linear operator from
$^{L\!}H^{\!\delta,\psi}(\Lambda)$ to $L^2_{\psi}(\Lambda)$.
	\end{enumerate}
(iv)-(vi) can be proved similarly.
\end{proof}

\begin{lemma} ($\psi$--fractional Poincar$\acute{e}$-Friedrichs inequalities)
\label{PFineq} The following two Poincar$\acute{e}$-Friedrichs-type inequalities
hold
\begin{equation*}
		\|v\|_{L^2_{\psi}(\Lambda)}\lesssim |v|_{^{L\!}H^{\!\delta,\psi}(\Lambda)},
		\ \ \ \forall v\in {}^{L\!}H^{\!\delta,\psi}(\Lambda).
	\end{equation*}
	\begin{equation*}
		\|v\|_{L^2_{\psi}(\Lambda)}\lesssim |v|_{^{R\!}H^{\!\delta,\psi}(\Lambda)},
	\ \ \ \forall v\in {}^{R\!}H^{\!\delta,\psi}(\Lambda).
	\end{equation*}
\end{lemma}
\begin{proof}
For all $v\in \,^{L\!}H^{\!\delta,\psi}(\Lambda),$ it follows from \eqref{deri_jifen}
and Lemma \ref{map_lem} that
	$$\|v\|_{L^2_{\psi}(\Lambda)}=\|I^{\!\delta,\psi}_tD^{\!\delta,\psi}_tv\|_{L^2_{\psi}(\Lambda)}\lesssim \|D^{\!\delta,\psi}_tv\|_{L^2_{\psi}(\Lambda)}=|v|_{^{L\!}H^{\!\delta,\psi}(\Lambda)}.$$
This proves the first inequality.
The second inequality can be derived similarly.
\end{proof}

One of the remarkable properties of the $\psi$--Riemann--Liouville fractional derivative is given in the following lemma.
\begin{lemma}\label{conjugate_lem1}
	For all $0<\delta<1,$ if $v\in H^{\!\delta,\psi}(\Lambda)$, $w\in C^{\infty}_0(\Lambda),$ then
	\begin{equation}\label{part_inte}
		\left( D_{t}^{\!\delta,\psi}v(t),~w(t)\right)_{L^2_{\psi}(\Lambda)}=\left( v(t),~_{t}D^{\!\delta,\psi}w(t)\right)_{L^2_{\psi}(\Lambda)}.
	\end{equation}
\end{lemma}

\begin{proof}
	By using integration by parts, we have
	\begin{align}
		\begin{split}
			&\left( D_{t}^{\!\delta,\psi}v(t),~w(t)\right)_{L^2_{\psi}(\Lambda)}\\
			=&\frac{1}{\Gamma(1-\delta)}\int_{\Lambda}\frac{1}{\psi^{\prime}(t)}
			\frac{d}{dt}\int^{t}_{a_\psi}\psi^{\prime}(s)(\psi(t)-\psi(s))^{-\delta}v(s)dsw(t)\psi^{\prime}(t)dt\\
			=&\frac{1}{\Gamma(1-\delta)}\big[w(t)\int^t_{a_\psi}\psi^{\prime}(s)(\psi(t)-\psi(s))^{-\delta}v(s)ds
			\big|_{\partial \Lambda}\\
			&-\int_{\Lambda}\int^{t}_{a_\psi}\psi^{\prime}(s)(\psi(t)-\psi(s))^{-\delta}v(s)dsdw(t)\big]\\
			=&\frac{-1}{\Gamma(1-\delta)}\int_{\Lambda}\int^{t}_{a_\psi}\psi^{\prime}(s)(\psi(t)-\psi(s))^{-\delta}v(s)dsw^{\prime}(t)dt\\
			=&\frac{-1}{\Gamma(1-\delta)}\int_{\Lambda}\int^{b_\psi}_s\left[\psi(t)-\psi(s)\right]^{-\delta}
			w^{\prime}(t)dtv(s)\psi^{\prime}(s)ds.
		\end{split}
	\end{align}
	Furthermore, a direct calculation gives
	\begin{align*}
		\begin{split}
			\frac{d}{ds}\int^{b_\psi}_s\psi^{\prime}(t)\left[\psi(t)-\psi(s)\right]^{-\delta}w(t)dt
			=&\frac{d}{ds}\big[w(t)\frac{\left[\psi(t)-\psi(s)\right]^{1-\delta}}{1-\delta}\big|^{b_\psi}_{s}\\
			&-\int^{b_\psi}_{s}\frac{\left[\psi(t)-\psi(s)\right]^{1-\delta}}{1-\delta}dw(t)\big]\\
			=&-\frac{d}{ds}\int^{b_\psi}_{s}\frac{\left[\psi(t)-\psi(s)\right]^{1-\delta}}{1-\delta}dw(t)\\
			=&\psi^{\prime}(s)\int^{b_\psi}_{s}\left[\psi(t)-\psi(s)\right]^{-\delta}w^{\prime}(t)dt.
		\end{split}
	\end{align*}
	Thus,
	\begin{align*}
		\begin{split}
			&\left( v(s),~_{s}D^{\!\delta,\psi}w(s)\right)_{L^2_{\psi}(\Lambda)}\\
			=&\frac{-1}{\Gamma(1-\delta)}\int_{\Lambda}\frac{1}{\psi^{\prime}(s)}
			\frac{d}{ds}\int^{b_\psi}_s\psi^{\prime}(t)\left[\psi(t)-\psi(s)\right]^{-\delta}w(t)dtv(s)\psi^{\prime}(s)ds\\
			=&\frac{-1}{\Gamma(1-\delta)}\int_{\Lambda}\frac{1}{\psi^{\prime}(s)}\psi^{\prime}(s)
			\int^{b_\psi}_s\left[\psi(t)-\psi(s)\right]^{-\delta}w^{\prime}(t)dtv(s)\psi^{\prime}(s)ds\\
			=&\frac{-1}{\Gamma(1-\delta)}\int_{\Lambda}\int^{b_\psi}_s\left[\psi(t)-\psi(s)\right]^{-\delta}
			w^{\prime}(t)dtv(s)\psi^{\prime}(s)ds\\
			=&\left( D_{t}^{\!\delta,\psi}v(t),~w(t)\right)_{L^2_{\psi}(\Lambda)}.
		\end{split}
	\end{align*}
This completes the proof.
\end{proof}
Based on a similar idea introduced in \cite{LX10}, the $\psi$--fractional derivative can be generalized as a distribution
to any $L^2_\psi(\Lambda)$ functions by using the integration by parts \eqref{part_inte}. That is,
for $v\in L_\psi^2(\Lambda)$,
the $\psi$--fractional derivative of $v$ in the distribution sense
is defined as the linear functional through
\bex\label{distribution}
\left\langle D_{t}^{\!\delta,\psi}v(t),~w(t)\right\rangle_{L^2_{\psi}(\Lambda)}=\left( v(t),~_{t}D^{\!\delta,\psi}w(t)\right)_{L^2_{\psi}(\Lambda)},\ \ \
\forall w\in C^{\infty}_0(\Lambda).
\eex
With this convention, we are able to derive,
by following the same lines as in \cite{LX10},
a key result which is crucial for the proof of well-posedness of the
variational problem. That is,
	for all $0<\delta<1,$ if $v, w\in H^{\!\delta/2,\psi}(\Lambda)$, then
\be\label{keyeq_con}
\left\langle  D_{t}^{\!\delta,\psi}v(t),~w(t)\right\rangle_{L^2_{\psi}(\Lambda)}=\left( D_{t}^{\!\delta/2,\psi}v(t),~_{t}D^{\!\delta/2,\psi}w(t)\right)_{L^2_{\psi}(\Lambda)}.
\ee

\begin{remark}
	It is worth noting that the $\psi$--fractional variational framework established in this section
	is valid for quite general function $\psi(t)$. The only assumption on $\psi(t)$ is its increasing differentiability and
$\psi^{\prime}(t)\neq 0$.
In what follows we will consider a special case $\psi(t)=t^{{1/\gamma}}$ to demonstrate how this
variational framework can be used to capture some singular solutions of fractional differential equations.
\end{remark}

\section{A spectral method for fractional ordinary differential equations}\label{Sect3}

As a simple application example, we consider
in this section the following initial value problem
\begin{equation}\label{Caputo_Problem}
	\left\{\begin{aligned}
		&^C\!D_s^\delta u(s)+\lambda u(s)= g(s), \ \lambda>0,\ s\in I,\\
		&u(a)=\phi.
	\end{aligned}\right.
\end{equation}
Here $0<\delta<1$, $^C\!D^\delta_s$
denotes the classical left--sided Caputo fractional operator defined in \eqref{classic_Caputo}.
This model problem frequently appears in the investigation of the TFDE \cite{Sak11,FX21}:
\bex
	^C\!D_s^\delta u(x,s)=\Delta u(x,s)+g(x,s),\ x\in \Omega,~s\in I,
\eex
where $\Omega$ is a spatial domain.
The solution of the TFDE can be expended in the space direction by using the eigenfunctions of the Laplacian operator $-\Delta$,
resulting in the equation \eqref{Caputo_Problem} with $\lambda$ being an eigenvalue of $-\Delta$.
It is seen that the model problem \eqref{Caputo_Problem} reflects the main difficulty of solving the TFDE, i.e., singularity
feature of the solution in the time direction.

Without loss of generality, we consider the homogeneous initial condition, i.e., $\phi\equiv 0$.
The case of non-homogeneous initial condition can be
handled by standard homogenization.  With $\phi\equiv 0$, the problem \eqref{Caputo_Problem} can be equivalently
written as \cite{YX21}
\be\label{RL_Problem}
\left\{\begin{aligned}
& D_s^\delta u(s)+\lambda u(s)= g(s), \ \lambda>0,\ s\in I,\\
& I^{1-\delta}_s u(a) = 0,
\end{aligned}\right.
\ee
where the operators $D_s^\delta$ and $I^{1-\delta}_s$ are defined in \eqref{lRD} and \eqref{lRI} with $\psi(s)=s$.

By the change of variable $s=\psi(t)$, and denoting
$v(t)=u(\psi(t)), f(t)=g(\psi(t))$,
the problem \eqref{RL_Problem} can be transformed into the following problem
\begin{equation}\label{problem-1}
	\left\{\begin{aligned}
		&\!D^{\!\delta,\psi}_tv(t)+\lambda v(t)= f(t), \ t\in \Lambda,\\
		&I^{1-\delta,\psi}_t v(a_\psi)=0.\end{aligned}\right.
\end{equation}

We propose and analyze below a spectral Galerkin method to solve the transformed
problem \eqref{problem-1} expressed in the weak form.
We first introduce the GJPs (see \cite{Guo2009Generalized, Shen2011Spectral}).
Define the shifted GJPs
\begin{equation}\label{Def_GJ}
{j}^{\alpha,-1}_n(t):=(1+x(t))J^{\alpha,1}_{n-1}(x(t)),~~\alpha>-1,~t\in \Lambda,~n=1,2,\dots,
\end{equation}
where $x(t)=\frac{2t-(a_\psi+b_\psi)}{b_\psi-a_\psi}$, $J^{\alpha,\beta}_{n}(x)$ are the classical $n$-th Jacobi polynomials, i.e., orthogonal
polynomials with respect to the weight function $\omega^{\alpha,\beta}(x):=(1-x)^\alpha(1+x)^\beta, \alpha,\beta>-1$,
$n=0,1,\dots$.

It can be checked that
\begin{equation}\label{Deri_GJ}
	\frac{d}{dt}{j}^{\alpha,-1}_n(t)=\frac{2n}{b_\psi-a_\psi}J^{\alpha+1,0}_{n-1}(x(t)),~~t\in \Lambda.
\end{equation}
Let $P_N$ be the standard polynomial space defined by
$$
P_N:=\mathrm{span}\{1,t,t^2,\ldots,t^N\}.
$$
Set the shifted polynomials space
\begin{equation*}
	V_N:=\mathrm{span}\{v\in P_N|v(0)=0\}=\mathrm{span}\left\{j^{\alpha,-1}_n(t),~~~t\in \Lambda,~n=1,\ldots,N\right\}.
\end{equation*}
Define the $L^2_{\omega^{\alpha,-1}}$--orthogonal projection $\pi^{\alpha,-1}_N$: $L^2_{\omega^{\alpha,-1}}(\Lambda)\rightarrow V_N$, such that for all $u\in L^2_{\omega^{\alpha,-1}}(\Lambda)$, $\pi^{\alpha,-1}_Nu\in V_N$ satisfies
\begin{equation*}
	\left(\pi^{\alpha,-1}_Nu-u,w\right)_{L^2_{\omega^{\alpha,-1}}(\Lambda)}=0,~~~\forall w\in V_N.
\end{equation*}
Define the non-uniform Jacobi-weighted Sobolev spaces as follows:
\bex
B^{m}_{\omega^{\alpha,-1}}(\Lambda)
:=\left\{v: \partial^{k}_{t}v\in L^{2}_{\omega^{\alpha+k,-1+k}}(\Lambda),0\leq k\leq m\right\}.
\eex
An approximation result of this projection operator is given in the following lemma.
\begin{lemma}\label{pro_err_lem}
~For any $u\in B^m_{\omega^{\alpha,-1}}(\Lambda), m\in \mathbb{N}, m\geq 1$,
	and $0\leq \mu\leq m,$ we have
	\begin{equation}\label{proj1}
		\left\|\pi^{\alpha,-1}_Nu-u\right\|_{B^{\mu}_{\omega^{\alpha,-1}}}\lesssim
		N^{\mu-m}\left\|\partial^m_tu\right\|_{L^2_{\omega^{\alpha+m,-1+m}}}.
	\end{equation}
\end{lemma}
\begin{proof}
This approximation result can be proved in the same way as for the projector $\pi^{-1,\beta}_N$ given
 in \cite{Guo2009Generalized}. We omit the details in order to limit the length of the paper. 	
\end{proof}	

The spectral approximation we propose for \eqref{problem-1} reads:
Find $v_N\in V_N$ such that
\begin{equation}\label{GSM}
	\mathcal{A}( v_N,~w_N)=F(w_N),~~~~\forall~w_N\in V_N,
\end{equation}
where
\bex
&&\mathcal{A}( v_N,~w_N)=\left(D^{\delta,\psi}_tv_N,~w_N\right)_{L^2_{\psi}(\Lambda)}
+\lambda\left( v_N,~w_N\right)_{L^2_{\psi}(\Lambda)},\\
&&F(w_N)=\left( f,~w_N\right)_{L^2_{\psi}(\Lambda)},
\eex
with $(\cdot, \cdot)_{L^2_{\psi}(\Lambda)}$ being defined in \eqref{innerp}.

\subsection{Well-posedness}

\begin{theorem}\label{Well-pose}
For any $f$ satisfying $I_t^{\delta/2,\psi} f\in L^2_{\psi}(\Lambda)$, $0<\delta<1$,
the spectral approximation problem \eqref{GSM} is well-posed.
Moreover, if $v_N$ is the solution of \eqref{GSM}, then it holds
	\begin{equation}\label{StabilityIneq}
		\|v_N\|_{H^{\delta/2,\psi}(\Lambda)}\lesssim \|I^{\delta/2,\psi}_t f\|_{L^2_{\psi}(\Lambda)}.
	\end{equation}
\end{theorem}
\begin{proof}
The proof makes use the classical Lax-Milgram Theorem, which consists in verifying the coercivity and
continuity of the bilinear form $\mathcal{A}(\cdot,\cdot)$.

Combining \eqref{keyeq_con} with the definition of
	$\left\|\cdot\right\|_{^{c\!}H^{\!\delta,\psi}(\Lambda)}$ gives:
	for all $v_N\in V_N$,
	\begin{align*}\label{Coe}
		\begin{split}
			\mathcal{A}\left(v_N,~v_N\right)
			=&\left(D_{t}^{\!\delta,\psi}v_N,~v_N\right)_{L^2_{\psi}(\Lambda)}
			+\lambda\left(v_N,~v_N\right)_{L^2_{\psi}(\Lambda)}\\
			=&\left(D_{t}^{\delta/2,\psi}v_N,
			~_{t}D^{\delta/2,\psi}v_N\right)_{L^2_{\psi}(\Lambda)}
			+\lambda\left(v_N,~v_N\right)_{L^2_{\psi}(\Lambda)}\\
			\gtrsim& \left|v_N\right|^2_{^{c\!}H^{\delta/2,\psi}(\Lambda)}.
		\end{split}
	\end{align*}
	Furthermore, the norm equivalence proved in Lemma~\ref{lem_equi} yields
	\begin{equation*}
		\mathcal{A}\left(v_N,~v_N\right)	
		\gtrsim \left|v_N\right|^2_{H^{\delta/2,\psi}(\Lambda)}.
	\end{equation*}
Then it follows from the fractional Poincar$\mathrm{\acute{e}}$-Friedrichs inequality in Lemma \ref{PFineq}:
	\begin{equation*}
		\mathcal{A}\left(v_N,~v_N\right)
		\gtrsim \left\|v_N\right\|^2_{H^{\delta/2,\psi}(\Lambda)}.
	\end{equation*}
By applying \eqref{keyeq_con} again, and using Cauchy-Schwarz inequality,
we obtain for all  $v_N, w_N\in V_N$,
	\begin{align*}
		\begin{split}
			\left|\mathcal{A}\left(v_N,~w_N\right)\right|
			\leq&\left|\left(D_{t}^{\!\delta,\psi}v_N,~w_N\right)_{L^2_{\psi}(\Lambda)}\right|
			+\lambda\left|\left(v_N,~w_N\right)_{L^2_{\psi}(\Lambda)}\right|\\
			\leq&\left|\left(D_{t}^{\delta/2,\psi}v_N,
			~_{t}D^{\delta/2,\psi}w_N\right)_{L^2_{\psi}(\Lambda)}\right|
			+\lambda\left|\left(v_N,~w_N\right)_{L^2_{\psi}(\Lambda)}\right|\\
			\leq& \left|v_N\right|_{^{L\!}H^{\delta/2,\psi}(\Lambda)}
			\left|w_N\right|_{^{R\!}H^{\delta/2,\psi}(\Lambda)}
			+\lambda\left\|v_N\right\|_{L^2_{\psi}(\Lambda)}
			\left\|w_N\right\|_{L^2_{\psi}(\Lambda)}.
		\end{split}
	\end{align*}
Finally, we derive from the norm equivalence and Lemma \ref{PFineq}:
	\begin{align*}
		\begin{split}
			\left|\mathcal{A}\left(v_N,~w_N\right)\right|
			\lesssim& \left|v_N\right|_{H^{\delta/2,\psi}(\Lambda)}
			\left|w_N\right|_{H^{\delta/2,\psi}(\Lambda)}\\
			\lesssim& \left\|v_N\right\|_{H^{\delta/2,\psi}(\Lambda)}
			\left\|w_N\right\|_{H^{\delta/2,\psi}(\Lambda)}.
		\end{split}
	\end{align*}
The well-posedness of \eqref{GSM} is thus proved.
	
The stability inequality \eqref{StabilityIneq} can be derived by using  \eqref{deri_jifen}, Lemma~\ref{lem_equi}, Lemma \ref{conjugate_lem1}, and the Cauchy-Schwarz inequality
as follows:
	\begin{align*}
		\begin{split}
			\left\|v_N\right\|^2_{H^{\delta/2,\psi}(\Lambda)}
			&\lesssim \mathcal{A}\left(v_N,~v_N\right)
			=\mathcal{F}(v_N)\\
			&= \left(D_t^{\delta/2,\psi}I_t^{\delta/2,\psi} f,~v_N\right)_{L^2_{\psi}(\Lambda)}
			= \left(I_t^{\delta/2,\psi} f,~_tD^{\delta/2,\psi}v_N\right)_{L^2_{\psi}(\Lambda)}\\
			&\leq \left\|I_t^{\delta/2,\psi}f\right\|_{L^2_{\psi}(\Lambda)}
			\|_tD^{\delta/2,\psi}v_N\|_{L^2_{\psi}(\Lambda)}\\
			&\leq \left\|I_t^{\delta/2,\psi}f\right\|_{L^2_{\psi}(\Lambda)}
			\left\|v_N\right\|_{H^{\delta/2,\psi}(\Lambda)}.
		\end{split}
	\end{align*}
This ends the proof of the theorem.
\end{proof}

\subsection{Error estimate}
In this subsection we present an error estimate for a specific transformation function, i.e.,
$\psi(t)=t^{1/{\gamma}}$.
Although this is the only case for which we derive the error estimate here, we are going to
see that this specific transformation can well smooth the time fractional diffusion equation, therefore
is a good fit for use of the spectral method.

Before carrying out the error analysis, we recall the following definition and lemma from \cite{Andersen1983Weighted}.
Define the integral operator $(P_{\delta})v(t)=\int^t_{-\infty}P_{\delta}(t,z)v(z)dz,$ where $P_{\delta}(t,z)\geq 0$ is not increasing in $t$ and not decreasing in $z$ in $\{(t,z)\in \mathbb{R}^2:z<t\}$. For two nonnegative functions
$\omega_1$ and $\omega_2$, we set
		\begin{align}\label{A_delta}
	A_{\delta, p, q}(t)=
	\left(\int^{\infty}_t\left[P_{\delta}(y,t)^{1/2}\omega_1(y)
	\right]^qdy\right)^{1/q}
	\left(\int^{t}_{-\infty}\left[P_{\delta}(t,y)^{-1/2}\omega_2(y)
	\right]^{-p^{\prime}}dy\right)^{1/{p^{\prime}}}.
\end{align}
\begin{lemma}\cite{Andersen1983Weighted}\label{key_lem}\
Let $\omega_1(t)\geq 0,$ $\omega_2(t)\geq 0,$
		$1/p+1/p^{\prime}=1.$
	If there exists a constant $C$ such that $A_{\delta, p, q}(t)\leq C, \forall t\in \mathbb{R}$. Then
	\begin{equation*}\label{AC}
		\left(\int^{\infty}_{-\infty}\left|\omega_1(t) (P_{\delta})v(t)\right|^qdt\right)^{1/q}
		\leq A_1C
		\left(\int^{\infty}_{-\infty}\left|\omega_2(t)v(t)\right|^pdt\right)^{1/p},
	\end{equation*}
	where
	$A_1=((p^{\prime}+q)/q)^{1/p^{\prime}}((p^{\prime}+q)/p^{\prime})^{1/q}$ if $1<p\leq q<\infty,$ and $A_1=1$ otherwise.
\end{lemma}

Next we prove two lemmas which are useful for the error estimation.

\begin{lemma}\label{key2_lem}
Let $P_{\delta}(t,~z)=\frac{1}{\Gamma(1-\delta)}(t^{1/\gamma}-z^{1/\gamma})^{-\delta}$,
 $\delta<\alpha<1, \Lambda_1=(0,1)$,
	\begin{align*}
		\begin{split}
			\omega_1(t)=
			\begin{cases}
				(1-t)^{\alpha/2}t^{1/{\gamma}-1/2},&t\in \Lambda_1,\\
				0,                                          &t\in \mathbb{R}\setminus\Lambda_1,
			\end{cases}~~
			\omega_2(t)=
			\begin{cases}
				(1-t)^{(1-\alpha)/2}, &t\in \Lambda_1,\\
				0,                                          &t\in \mathbb{R}\setminus\Lambda_1.
			\end{cases}
		\end{split}
	\end{align*}
For any differentiable function  $v(t)$ defined in $\Lambda_1$, it holds
	\begin{equation*}
		\int_{\Lambda_1}\left|\omega_1(t)( P_{\delta}v^{\prime})(t)\right|^2dt\lesssim
		\int_{\Lambda_1}\left|\omega_2(t)(v^{\prime})(t)\right|^2dt.
	\end{equation*}
\end{lemma}

\begin{proof}
First, it is easy to check that $P_{\delta}(t,z)\geq 0$ is not increasing in $t$ and not decreasing in $z$.
We extend $v(t)$ with zero outside of $\Lambda_1$.
Taking $p=q=2$ in \eqref{A_delta} gives
	\begin{align*}
		\begin{split}
			&A_{\delta,2,2}(t)\\
			=&\left(\int^1_t\left[(P_{\delta}(y,t))^{1/2}\omega_1(y)\right]^2dy\right)^{1/2}
			\left(\int^t_0\left[(P_{\delta}(t,y))^{-1/2}\omega_2(y)\right]^{-2}dy\right)^{1/2}\\
			=&\left(\int^1_t\big(y^{1/\gamma}-t^{1/\gamma}\big)^{-\delta}(1-y)^{\alpha}y^{2/\gamma-1}dy\right)^{1/2}
			\left(\int^t_0\big(t^{1/\gamma}-y^{1/\gamma}\big)^{-\delta}(1-y)^{\alpha-1)}dy\right)^{1/2}.
		\end{split}
	\end{align*}
	A direct calculation shows
\begin{align*}
	\begin{split}
		&\int^1_t\left(y^{1/\gamma}-t^{1/\gamma}\right)^{-\delta}(1-y)^{\alpha}y^{2/\gamma-1}dy\\
		\leq&\int^1_t\left(y^{1/\gamma}-t^{1/\gamma}\right)^{-\delta}y^{2/\gamma-1}dy\\
		=&\gamma\int^1_t\left(y-t\right)^{-\delta}ydy\\
		\leq&\frac{\gamma}{1-\delta}.
	\end{split}
\end{align*}
Notice $1-y^{\gamma}\geq \gamma(1-y)$ for $\gamma\in (0,1], y\in (0,1),$ we have
\begin{align*}
	\begin{split}
		&\int^t_0\left(t^{1/\gamma}-y^{1/\gamma}\right)^{-\delta}(1-y)^{\alpha-1}dy\\
		=&\gamma\int^t_0(t-y)^{-\delta}(1-y^{\gamma})^{\alpha-1}y^{\gamma-1}dy\\
		\leq&\gamma^{\alpha}\int^t_0(t-y)^{-\delta}(1-y)^{\alpha-1}y^{\gamma-1}dy\\
		\leq&\gamma^{\alpha}B(\alpha-\delta,\gamma)t^{\alpha-\delta+\gamma-1}.
	\end{split}
\end{align*}
Combining the above estimates yields
	\begin{equation*}
		A_{\delta,2,2}(t)\leq C<\infty.
	\end{equation*}
Then we conclude by using Lemma \ref{key_lem}.
\end{proof}

\begin{lemma}\label{key3_lem}
	Assume $v\in\ _0H^{\delta/2,\psi}(\Lambda_1)\cap B^m_{\omega^{-\alpha,-1}}(\Lambda_1),$ where $0<\delta<1,$ $\psi(t)=t^{1/{\gamma}},$ $m\geq 1,$  $\delta<\alpha<1$. Then we have
	\begin{enumerate}
		\renewcommand\theenumi{\roman{enumi}}
		\renewcommand{\labelenumi}{$\mathrm{\left(\theenumi\right)}$}
		\item $\left\|v\right\|_{L^2_{\psi}(\Lambda_1)}\lesssim \left\|v\right\|_{L^2_{\omega^{-\alpha,-1}}(\Lambda_1)}$,
		\item $\left(^C\!D^{\!\delta,\psi}_tv,~v\right)_{L^2_{\psi}(\Lambda_1)}\lesssim
		\left\|v^{\prime}\right\|_{L^2_{\omega^{1-\alpha,0}}(\Lambda_1)}
		\left\|v\right\|_{L^2_{\omega^{-\alpha,-1}}(\Lambda_1)}$.
	\end{enumerate}
\end{lemma}
\begin{proof}
	\begin{enumerate}
		\renewcommand\theenumi{\roman{enumi}}
		\renewcommand{\labelenumi}{$\mathrm{\left(\theenumi\right)}$}
		\item Noticing ${1/\gamma}t^{1/{\gamma}-1}\leq (1-t)^{-\alpha}t^{-1},~\forall t\in \Lambda_1,$ we have
		$$\left\|v\right\|^2_{L^2_{\psi}(\Lambda_1)}=
		{1/\gamma}\int_{\Lambda_1}v^2t^{1/{\gamma}-1}dt
		\leq \int_{\Lambda_1}v^2(1-t)^{-\alpha}t^{-1}dt
		=\left\|v\right\|_{L^2_{\omega^{-\alpha,-1}}(\Lambda_1)}.$$
		\item
		By a direct computation, we get
		\begin{align*}
			\begin{split}
				\left(^C\!D^{\!\delta,\psi}_tv,~v\right)_{L^2_{\psi}(\Lambda_1)}
				=&\frac{1}{\Gamma(1-\delta)}\int_{\Lambda_1}\int^t_0\left(t^{1/{\gamma}}-z^{1/{\gamma}}\right)
				^{-\delta}v^{\prime}(z)dzv(t)\psi^{\prime}(t)dt\\
				=&{1/\gamma}\int_{\Lambda_1}(P_{\delta}v^{\prime})(t)\omega_1(t)v(t)(1-t)^{-\alpha/2}t^{-1/2}dt\\
				\leq&{1/\gamma}\left[\int_{\Lambda_1}\left|\omega_1(t)(P_{\delta}v^{\prime})(t)\right|^2dt\right]^{1/2}
				\left\|v\right\|_{L^2_{\omega^{-\alpha,-1}}(\Lambda_1)}\\
				\lesssim&\left(\int_{\Lambda_1}\left|\omega_2(t)v^{\prime}(t)\right|^2dt\right)^{1/2}
				\left\|v\right\|_{L^2_{\omega^{-\alpha,-1}}(\Lambda_1)}\\
				=&\left\|v^{\prime}\right\|_{L^2_{\omega^{1-\alpha,0}}(\Lambda_1)}
				\left\|v\right\|_{L^2_{\omega^{-\alpha,-1}}\left(\Lambda_1\right)}.
			\end{split}
		\end{align*}
		Here, the Cauchy--Schwarz inequality and Lemma~\ref{key2_lem} have been used.
	\end{enumerate}
\end{proof}

With the above preparation, we are now in a position to derive the error estimate.
\begin{theorem}
Let $v$ be the solution of problem \eqref{problem-1}, $v_N$ is the solution of the problem \eqref{GSM}.
Suppose $v\in\ _0H^{\delta/2,\psi}(\Lambda)\cap B^m_{\omega^{-\alpha,-1}}(\Lambda),$ where $0<\delta<1,$ $\psi(t)=t^{1/{\gamma}},$ $m\geq 1,$  $\delta<\alpha<1$. Then the following error estimate holds:
	\begin{equation}\label{err_esti}
		\|v-v_N\|_{H^{\delta/2,\psi}(\Lambda)}\lesssim N^{1/2-m}\left\|\partial^m_tv\right\|_{L^2_{\omega^{-\alpha+m,-1+m}}(\Lambda)}+N^{-m}\|\partial^m_tv\|_{L^2_{\omega^{-\alpha+m,-1+m}}(\Lambda)}.
	\end{equation}
\end{theorem}
\begin{proof}
	It follows from \eqref{problem-1}, \eqref{GSM}, and C\'ea lemma that
	\begin{equation*}
		\|v-v_N\|_{H^{\delta/2,\psi}(\Lambda)}\leq \mathop{\mathrm{inf}}\limits_{w_N\in V_N}\|v-w_N\|_{H^{\delta/2,\psi}(\Lambda)}
		\leq\|v-\pi_N^{-\alpha,-1}v\|_{H^{\delta/2,\psi}(\Lambda)}.
	\end{equation*}
Furthermore it is not difficult to derive
\be\label{above1}
\|v-\pi_N^{-\alpha,-1}v\|_{H^{\delta/2,\psi}(\Lambda)}\lesssim \|v-\pi_N^{-\alpha,-1}v\|_{H^{\delta/2,\psi}(\Lambda_1)}.
\ee
Then, it follows from the definition of the norm $\|\cdot\|_{H^{\delta/2,\psi}(\Lambda)},$  the equivalence of $\psi$--fractional norms,
 and the relationship \eqref{C-RL}:
	\begin{align*}
	\begin{split}
		&\|v-\pi_N^{-\alpha,-1}v\|_{H^{\delta/2,\psi}(\Lambda_1)}\\\lesssim&
		\left|\left(D^{\!\delta,\psi}_t(v-\pi_N^{-\alpha,-1}v),~v-\pi_N^{-\alpha,-1}v\right)
		_{L^2_{\psi}(\Lambda_1)}\right|^{1/2}+\left\|v-\pi_N^{-\alpha,-1}v\right\|_{L^2_{\psi}(\Lambda_1)}\\=& \left|\left(^C\!D^{\!\delta,\psi}_t(v-\pi_N^{-\alpha,-1}v),~v-\pi_N^{-\alpha,-1}v\right)
		_{L^2_{\psi}(\Lambda_1)}\right|^{1/2}+\left\|v-\pi_N^{-\alpha,-1}v\right\|_{L^2_{\psi}(\Lambda_1)}.
	\end{split}
\end{align*}
Using Lemma \ref{key3_lem} gives
	\begin{align}\label{above3}
		\begin{split}
			&\|v-\pi_N^{-\alpha,-1}v\|_{H^{\delta/2,\psi}(\Lambda_1)}\\
			\lesssim& \left\|(v-\pi_N^{-\alpha,-1}v)'\right\|^{1/2}_{L^2_{\omega^{1-\alpha,0}}(\Lambda_1)}
			\left\|v-\pi_N^{-\alpha,-1}v\right\|^{1/2}_{L^2_{\omega^{-\alpha,-1}}\left(\Lambda_1\right)}+\left\|v-\pi_N^{-\alpha,-1}v\right\|_{L^2_{\omega^{-\alpha,-1}}\left(\Lambda_1\right)}.
		\end{split}
	\end{align}
	Finally, the desired estimate follows from combining \eqref{above1}, \eqref{above3}, and Lemma \ref{pro_err_lem}.
\end{proof}
\subsection{Implementation}
We discuss the implementation issue of the spectral approximation \eqref{GSM}.
The key is find efficient way to form the stiffness matrix $\textbf{S}$, those entries are
$$\textbf{S}_{mn}:=\left(\!D^{\delta,\psi}_tj^{\alpha,-1}_{n}(t),~j^{\alpha,-1}_{m}(t)\right)_{L^2_{\psi}(\Lambda)}
=\left(^C\!D^{\delta,\psi}_tj^{\alpha,-1}_{n}(t),~j^{\alpha,-1}_{m}(t)\right)_{L^2_{\psi}(\Lambda)}$$
for $m,~n=1,2,\cdots,N.$
We compute the entries $\textbf{S}_{mn}$ by using \eqref{Def_GJ} and \eqref{Deri_GJ} as follows:
\begin{align}\label{reStiff}
	\begin{split}
		&\left(^C\!D^{\delta,\psi}_tj^{\alpha,-1}_{n}(t),~j^{\alpha,-1}_{m}(t)\right)_{L^2_{\psi}(\Lambda)}\\
		=&\left(\frac{1}{\Gamma(1-\delta)}\int_0^s\left(s^{1/\gamma} -z^{1/\gamma}\right)^{-\delta}\frac{d}{dz}j^{\alpha,-1}_{n}(z)dz,~j^{\alpha,-1}_{m}(s)\right)_{L^2_{\psi}(\Lambda)}\\
		=&\frac{4nT^{1-\delta}}{\gamma\Gamma(1-\delta)}\!\!\int^1_0\!\!\!\int_0^1\!\!\!(1-\tau)^{-\delta}\big(\frac{1 -\tau^{{1/\gamma}}}{1-\tau}\big)^{-\delta}J^{\alpha+1,0}_{n}(2\tau s-1)d\tau
		\cdot s^{(1-\delta)/\gamma+1}J^{\alpha,1}_{m}(2s-1)ds\\
		=&\frac{4nT^{1-\delta}}{\gamma\Gamma(1-\delta)}\sum^{N}_{i=0}\sum^{N}_{j=0}
		\left(\frac{1-\hat{\eta}_j^{1/\gamma}}{1-\hat{\eta}_j}\right)^{-\delta}
		J^{\alpha+1,0}_n(2\eta_i\hat{\eta}_j-1)\hat{\omega}_jJ^{\alpha,1}_m(2\eta_i-1)\omega_i,
	\end{split}
\end{align}
where the Gauss quadrature point sets $\{\eta_i\}^{N}_{i=0}$ and $\{\hat{\eta}_j\}^{N}_{j=0}$ are zeros of the shifted Jacobi polynomials $J^{0,~(1-\delta)/\gamma+1}_{N}(x(t)),$ $J^{-\delta,~0}_{N}(x(t)),$ respectively, and
$\{\omega_i\}^{N}_{i=0},$ $\{\hat{\omega}_j\}^{N}_{j=0}$ are the associated weights.
Note that in our calculation, $\gamma$ is set to be $1/r,$ with $r$ being the positive integers
so that $\frac{1-\tau^{1/\gamma}}{1-\tau}=\sum^{r-1}_{k=0}\tau^k,$ which hasn't singularity. The singular parts
$(1-\tau)^{-\delta}$ and $s^{(1-\delta)/\gamma+1}$ do not appear in the numerical quadrature since they are treated as
the associated weights of the Jacobi polynomials.
Denote
\begin{equation*}
	f_m=(f,j^{\alpha,-1}_{m}(x)),~~~\mathbf{f}=(f_1,\ldots,f_N)^{T};
\end{equation*}
\begin{equation*}
	v_N=\sum^N_{n=1}\tilde{v}_nj^{\alpha,-1}_{n}(x),~~~\mathbf{v}=
	(\tilde{v}_1,\ldots,\tilde{v}_N)^{T};
\end{equation*}
\begin{equation*}
	\textbf{M}=(\textbf{M}_{mn})_{1\leq m,n\leq N} \mbox{ with } \textbf{M}_{mn}=\left(j^{\alpha,-1}_{n}(t),~j^{\alpha,-1}_{m}(t)\right)_{L^2_{\psi}(\Lambda)}.
\end{equation*}
Then the matrix form of the problem \eqref{GSM} reads:
\begin{equation*}
	(\textbf{S}+\lambda \textbf{M})\mathbf{v}=\mathbf{f}.
\end{equation*}

\begin{remark}
In \cite{Hou2018Muntz}, the authors proposed a M\"untz spectral method based on the
M\"untz polynomial space
$\mathrm{span}\{1,t^{\gamma},t^{2\gamma},\ldots,t^{N\gamma}\}$ for the fractional differential
equation.
It can be verified that, with the particular choice of the transformation $\psi(t)=t^{1/\gamma}$,
the current method is equivalent to the one in \cite{Hou2018Muntz} in the
sense that the solution $u_N(s)$ computed from the M\"untz spectral method is linked to
the solution $v_N(t)$ of the $\psi$--spectral method through $v_N(t)=u_N(\psi(t))$.
However, it is worth to note that the numerical analysis of the two methods
was conducted using two quite different frameworks.
The new approach in the current work not only provides an alternative tool for numerical analysis of the M\"untz spectral methods proposed
in \cite{Hou2017A,Hou2018Muntz},
but also provide a guideline for the selection of parameter $\gamma$.
The main goal is to choose a suitable transformation function $\psi$ such that $v(\cdot)=u(\psi(\cdot))$
is as smooth as possible.
\end{remark}

\section{ Application to the time fractional subdiffusion equations}\label{Sect4}
Let $\Omega:=(-1,1)^d,~d=1,2,3$. Consider the following time fractional diffusion equation:
\begin{equation}\label{TFDE}
	^C\!D_s^\delta u(x,s)=\Delta u(x,s)-u(x,s)+g(x,s),\ x\in \Omega,~s\in I
\end{equation}
subject to the initial and boundary conditions
\begin{equation}\label{Ini}
	u(x,0)=0,\ x\in \Omega,
\end{equation}
\begin{equation}\label{Bound}
	u(x,s)=0,\ x\in \partial\Omega,~s\in I.
\end{equation}
We obtain the following transformed equation by applying the
transformation $s=\psi(t)$ in the time direction:
\begin{equation}\label{psi-TFDE}
D_t^{\delta,\psi} v(x,t)=\Delta v(x,t)-v(x,t)+f(x,t),\ x\in \Omega,~t\in \Lambda.
\end{equation}

For the Sobolev space $X$ with norm $\|\cdot\|_{X},$ let
\begin{align*}
	\begin{split}
		H^{\!\delta,\psi}(\Lambda,X):=\left\{v;\|v(\cdot,t)\|_X\in H^{\!\delta,\psi}(\Lambda)\right\},\ \delta\geq 0
	\end{split}
\end{align*}
endowed with the norm
$$\|v\|_{H^{\!\delta,\psi}(\Lambda,X)}:=\left\|\left\|v(\cdot,t)\right\|_X\right\|_{\delta,\Lambda}.$$
Let $\mathcal{O}=\Lambda\times \Omega,$
\bex
\mathcal{H}^{\!\delta,\psi}(\mathcal{O}):={}_0H^{\!\delta,\psi}(\Lambda,L^2(\Omega))\cap L^2_{\psi}(\Lambda,H^1_0(\Omega)),
\\
\mathcal{B}^m_{\omega^{\alpha,\beta}}\left(\mathcal{O}\right):=B^m_{\omega^{\alpha,\beta}}(\Lambda,L^2(\Omega))
\cap L^2_{\omega^{\alpha,\beta}}(\Lambda,H^1_0(\Omega)),
\eex
equipped respectively with the norms
\bex
\|v\|_{\mathcal{H}^{\!\delta,\psi}(\mathcal{O})}:=\Big(\left\|v\right\|^2_{H^{\!\delta,\psi}(\Lambda,L^2(\Omega))}+\left\|v\right\|^2_{L^2_{\psi}(\Lambda,H^1_0(\Omega))}\Big)^{1/2},\\
\left\|v\right\|_{\mathcal{B}^m_{\omega^{\alpha,\beta}}(\mathcal{O})}:=\Big(\left\|v\right\|^2_{B^m_{\omega^{\alpha,\beta}}(\Lambda,L^2(\Omega))}+
\left\|v\right\|^2_{L^2_{\omega^{\alpha,\beta}}(\Lambda,H^1_0(\Omega))}\Big)^{1/2}.
\eex

Consider the following variational formulation of \eqref{psi-TFDE}:
\begin{equation}\label{weak_formula_1}
	\mathcal{A}(v,w)=\mathcal{F}(w),~~~\forall w\in \mathcal{H}^{\delta/2,\psi}(\mathcal{O}),
\end{equation}
where
\bex
\mathcal{A}(v,w):=(D_{t}^{\!\delta,\psi}v,w)_{\mathcal{O}}
+(\nabla v, \nabla w)_{\mathcal{O}}+(v,w)_{\mathcal{O}}, \
\mathcal{F}(w):=(f,w)_{\mathcal{O}}
\eex
with $(v,w)_{\mathcal{O}}:=\int_{\Lambda}\int_{\Omega}vw\psi^{\prime}(t) dxdt$.

Similar to Theorem \ref{Well-pose}, we can establish the coercivity and continuity of the bilinear
form $\mathcal{A}(v,w)$ in the space
$\mathcal{H}^{\delta/2,\psi}(\mathcal{O})\times \mathcal{H}^{\delta/2,\psi}(\mathcal{O})$, and therefore
the well-posedness of the weak problem \eqref{weak_formula_1} for any
$f\in \mathcal{H}^{\delta/2,\psi}(\mathcal{O})^{\prime}$ (the dual space of $\mathcal{H}^{\delta/2,\psi}(\mathcal{O})$), together with the stability estimate
\bex
\left\|v\right\|_{\mathcal{H}^{\delta/2,\psi}(\mathcal{O})}\lesssim \left\|f\right\|_{\mathcal{H}^{\delta/2,\psi}(\mathcal{O})^{\prime}}.
\eex

We now propose a space-time Galerkin spectral method to discretize \eqref{weak_formula_1}. For the time variable, we follow the approach of the previous section. For the space variable, we use standard Legendre polynomials. Let
\begin{align}
	\begin{split}
		\phi_k(x)=c_k(L_k(x)-L_{k+2}(x)),\ \ c_k=\frac{1}{\sqrt{4k+6}},\\
		a_{jk}=(\partial_x\phi_k(x),\partial_x\phi_j(x)),\ b_{jk}=(\phi_k(x),\phi_j(x)),
	\end{split}
\end{align}
where $L_n(x)$ is the $n$-th degree Legendre polynomial.
Then \cite{Shen1994Efficient}
\begin{align}
	\begin{split}
		a_{jk}=
		\begin{cases}
			1,\ &k=j\\
			0,\ &k\neq j
		\end{cases}
		,\
		b_{jk}=b_{kj}=
		\begin{cases}
			c_kc_j(\frac{2}{2j+1}+\frac{2}{2j+5}),\ &k=j\\
			-c_kc_j\frac{2}{2k+1},\ &k=j+2\\
			0,\ &\mathrm{otherwise}
		\end{cases}
		,
	\end{split}
\end{align}
Set the polynomial space
\begin{equation*}
	\mathcal{P}_M=\mathrm{span}\{\phi_0(x),\ \phi_1(x),\ \ldots,\ \phi_{M-2}(x)\}.
\end{equation*}
The space-time Galerkin spectral method for \eqref{weak_formula_1} is to seek $v_L(x,t)\in \mathcal{P}_M\otimes V_N$ such that
\begin{equation}\label{CalDouLi}
	\mathcal{A}(v_L,w)=\mathcal{F}(w),~~~\forall~w\in \mathcal{P}_M\otimes V_{N}.
\end{equation}

The error estimate is given in the following theorem without proof.
\begin{theorem}
Let $v$ be the solution of problem \eqref{weak_formula_1}, $v_L$ is the solution of problem \eqref{CalDouLi}.
Suppose $v\in \mathcal{H}^{\delta/2,\psi}(\mathcal{O})\cap B^m_{\omega^{-\alpha,-1}}(\Lambda,H^{\sigma}(\Omega)),$ where $0<\delta<1,$ $\psi(t)=t^{1/{\gamma}},$ $m\geq 1,$  $\delta<\alpha<1$. Then the following error estimate holds:
	\begin{align*}
		\begin{split}
			\|v-v_N\|_{\mathcal{H}^{\delta/2,\psi}(\mathcal{O})}\lesssim&  N^{1/2-m}\big\|\partial^m_tv\big\|_{L^2_{\omega^{-\alpha+m,-1+m}}(\Lambda,L^2(\Omega))}\\
			&+N^{1/2-m}M^{-\sigma}\left\|\partial^m_tv\right\|_{L^2_{\omega^{-\alpha+m,-1+m}}(\Lambda,H^{\sigma}(\Omega))}\\
			&+M^{-\sigma}\left(\left\|D^{\delta/2,\psi}_tv\right\|_{L^2_{\omega^{-\alpha,-1}}(\Lambda,H^{\sigma}(\Omega))}+\big\|v\big\|_{L^2_{\omega^{-\alpha,-1}}(\Lambda,H^{\sigma}(\Omega))}\right)\\
			&+M^{1-\sigma}\big\|v\big\|_{L^2_{\omega^{-\alpha,-1}}(\Lambda,H^{\sigma}(\Omega))}+N^{-m}\left\|\partial^m_tv\right\|_{L^2_{\omega^{-\alpha+m,-1+m}}(\Lambda,H^1_0(\Omega))}.
		\end{split}
	\end{align*}
\end{theorem}

In the implementation,
we express the solution $v_L$ of Eq. \eqref{CalDouLi} using the modal basis as follows:
\bex
	v_L(x,t)=\sum^{M-2}_{m=0}\sum^N_{n=1}\hat{v}_{nm}\phi_m(x)j^{\alpha,-1}_{n}(t).
\eex
Substituting this expression into \eqref{CalDouLi}, and taking $w=\phi_p(x)j^{\alpha,-1}_{q}(t)$, we obtain
\begin{align*}
	\begin{split}
		\sum^{M-2}_{m=0}\sum^N_{n=1}\hat{v}_{nm}\Big\{\left(\phi_m,\phi_p\right)_{L^2(\Omega)}\left(D^{\!\delta,\psi}_tj^{\alpha,-1}_{n},j^{\alpha,-1}_{q}\right)_{L^2_\psi(\Lambda)}+\left(\phi^{\prime}_m,\phi^{\prime}_p\right)_{L^2(\Omega)}\left(j^{\alpha,-1}_{n},j^{\alpha,-1}_{q}\right)_{L^2_\psi(\Lambda)}\\+\left(\phi_m,\phi_p\right)_{L^2(\Omega)}
		\left(j^{\alpha,-1}_{n},j^{\alpha,-1}_{q}\right)_{L^2_\psi(\Lambda)}\Big\}
		=\left(f,\phi_pj^{\alpha,-1}_{q}\right)_{\mathcal{O}}.
	\end{split}
\end{align*}
Denote
\begin{align*}
	\begin{split}
		&\mathbf{B}=(b_{jk})_{0\leq j,k\leq M-2},\\
		&f_{nm}=\left(f,\phi_m(x)j^{\alpha,-1}_{n}\left(t\right)\right)_{\mathcal{O}},\ \ \ \textbf{F}=(f_{nm})_{1\leq n\leq N,\ 0\leq m\leq M-2},\\
		&\mathbf{S}=\left(\mathbf{S}_{nq}\right)_{1\leq n,q\leq N},\ \mathbf{M}=\left(\mathbf{M}_{nq}\right)_{1\leq n,q\leq N},\ \textbf{V}=\left(\hat{v}_{nm}\right)_{1\leq n\leq N,\ 0\leq m\leq M-2}.
	\end{split}
\end{align*}
Using the above notations, \eqref{CalDouLi} can be written under the following matrix form:
\bex
	\mathbf{S}\textbf{V}\mathbf{B}+\mathbf{M}\textbf{V}+\mathbf{M}\textbf{V}\mathbf{B} = \textbf{F}.
\eex

\section{Numerical examples}\label{Sect5}
In this section, we present some numerical examples to illustrate the high accuracy of the proposed method
based on GJPs in solving problem \eqref{Caputo_Problem} with smooth and nonsmooth solutions. In particular, we test the accuracy of the proposed method when the exact solution is unknown.
The space-time spectral method based on GJPs and Legendre polynomials presented in Sect. \ref{Sect4}
will also be tested for the two-dimensional time fractional subdiffusion equation.
The time interval is set to $[a_\psi,b_\psi]:=[0,2^\gamma]$. Note that $\psi(t)=t^{1/\gamma}~(0<\gamma\leq 1)$ in the following examples.
\begin{example}(Smooth solution)\label{exam1}
	In this test, we choose the fabricated exact solution $u(s)=s^{2}$. Naturally, in this case, we take $\gamma=1.$
\end{example}
The main purpose of this example is to check the high accuracy of the proposed Galerkin spectral scheme
\eqref{GSM} for smooth solutions. The computed results are presented in Table \ref{ErrTab1},
from which we observe that the numerical solutions for some different $\delta$ reach the machine accuracy
with small polynomial degree $N$.

\begin{table}[h]
	\centering
	\caption{(Example 1) $L^{\infty}$- and $L^2$-errors versus $N$ and different $\delta$.}
	\begin{tabular}{|c|c|c|c|c|c|c|}
		\hline
		& \multicolumn{2}{c|}{$\delta=0.1$} & \multicolumn{2}{c|}{$\delta=0.5$} & \multicolumn{2}{c|}{$\delta=0.9$} \\
		\hline
		$N$ & $\|v-v_N\|_{L^{\infty}}$ & $\|v-v_N\|_0$ & $\|v-v_N\|_{L^{\infty}}$ & $\|v-v_N\|_0$ & $\|v-v_N\|_{L^{\infty}}$ &    $\|v-v_N\|_0$ \\
		\hline
		2 & 3.3307e-16  & 4.0030e-16 & 1.1102e-16 & 1.3878e-16 & 1.4433e-15 & 1.4647e-15 \\
		4 & 1.3323e-15 & 1.4767e-15 & 1.5543e-15 & 1.6812e-15 & 1.1990e-14 & 1.4989e-14 \\
		\hline
	\end{tabular}\label{ErrTab1}
\end{table}
\begin{example}(Nonsmooth solution)\label{exam2}
Consider problem \eqref{Caputo_Problem} with the fabricated exact solution $u(s)=s^{\sigma}$
for two values of $\sigma: 3/5, \sqrt{2}/2$.
\end{example}

\begin{figure}[htbp]
	\centering
	\includegraphics[scale=0.35]{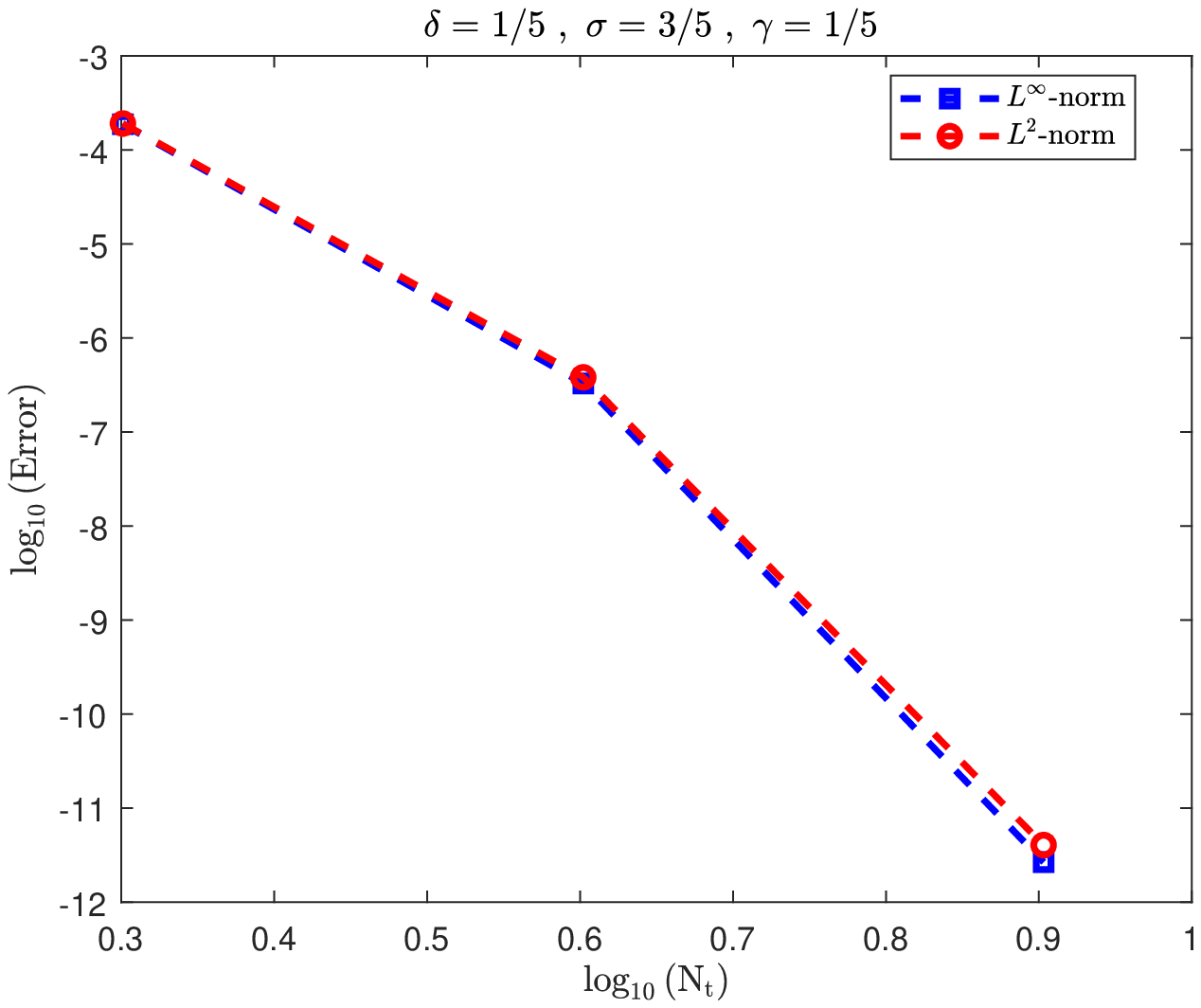}
	\includegraphics[scale=0.35]{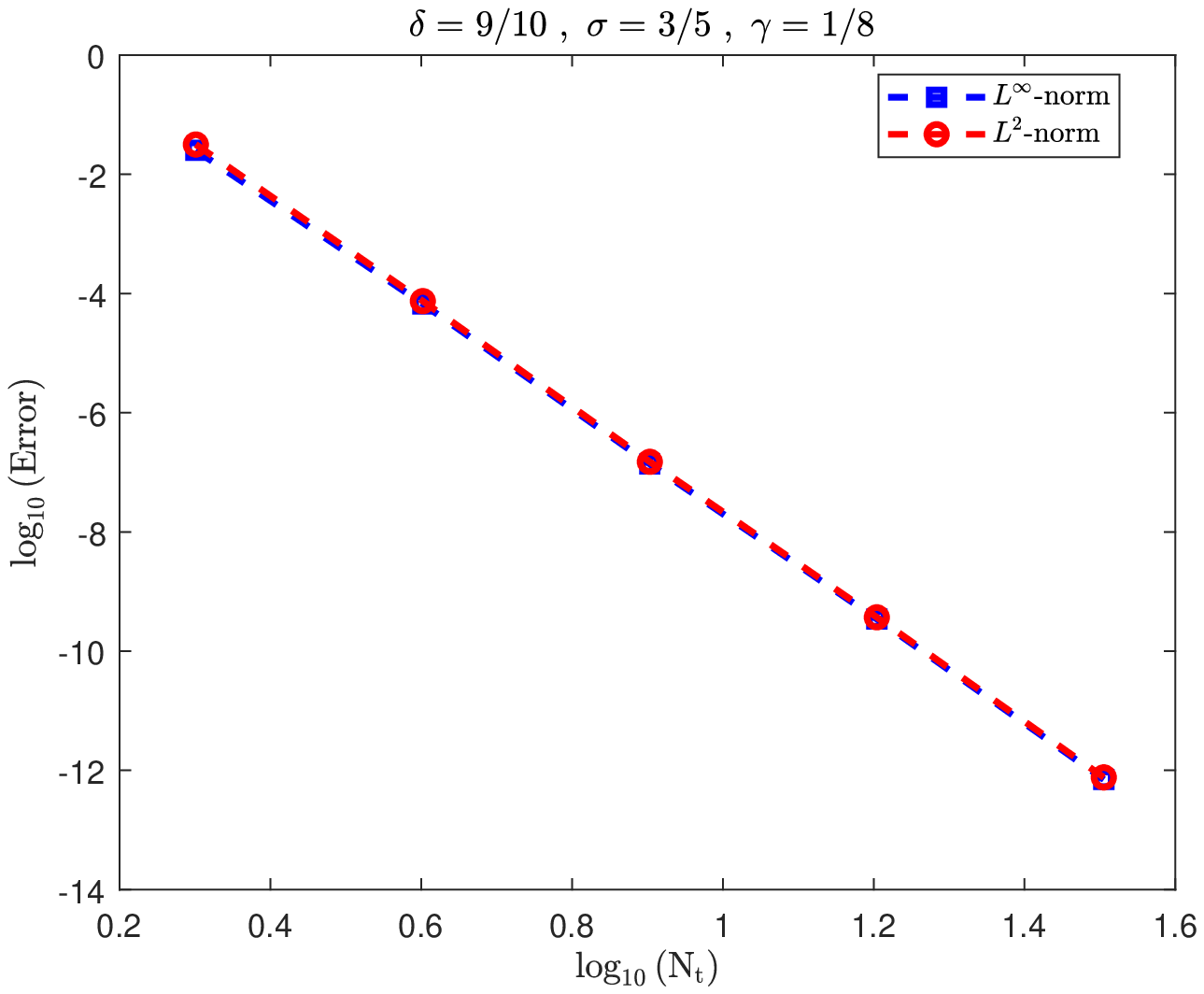}\\
	\caption{(Example 2) $L^{\infty}$- and $L^2$-errors in log scale versus the polynomial degree for fractional power solution.}\label{ErrPlot1}
\end{figure}

\begin{figure}[htbp]
	\centering
	\includegraphics[scale=0.35]{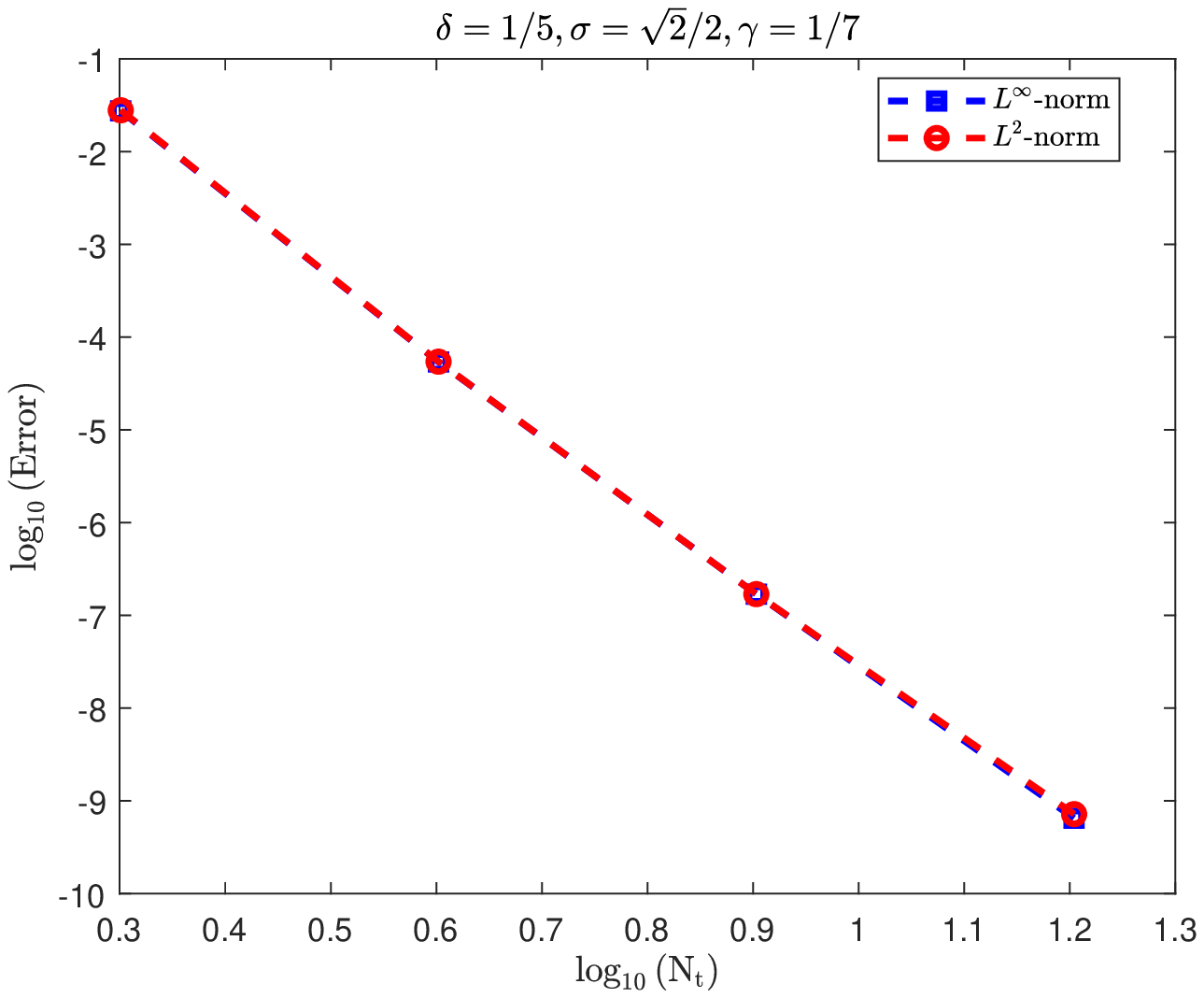}
	\includegraphics[scale=0.35]{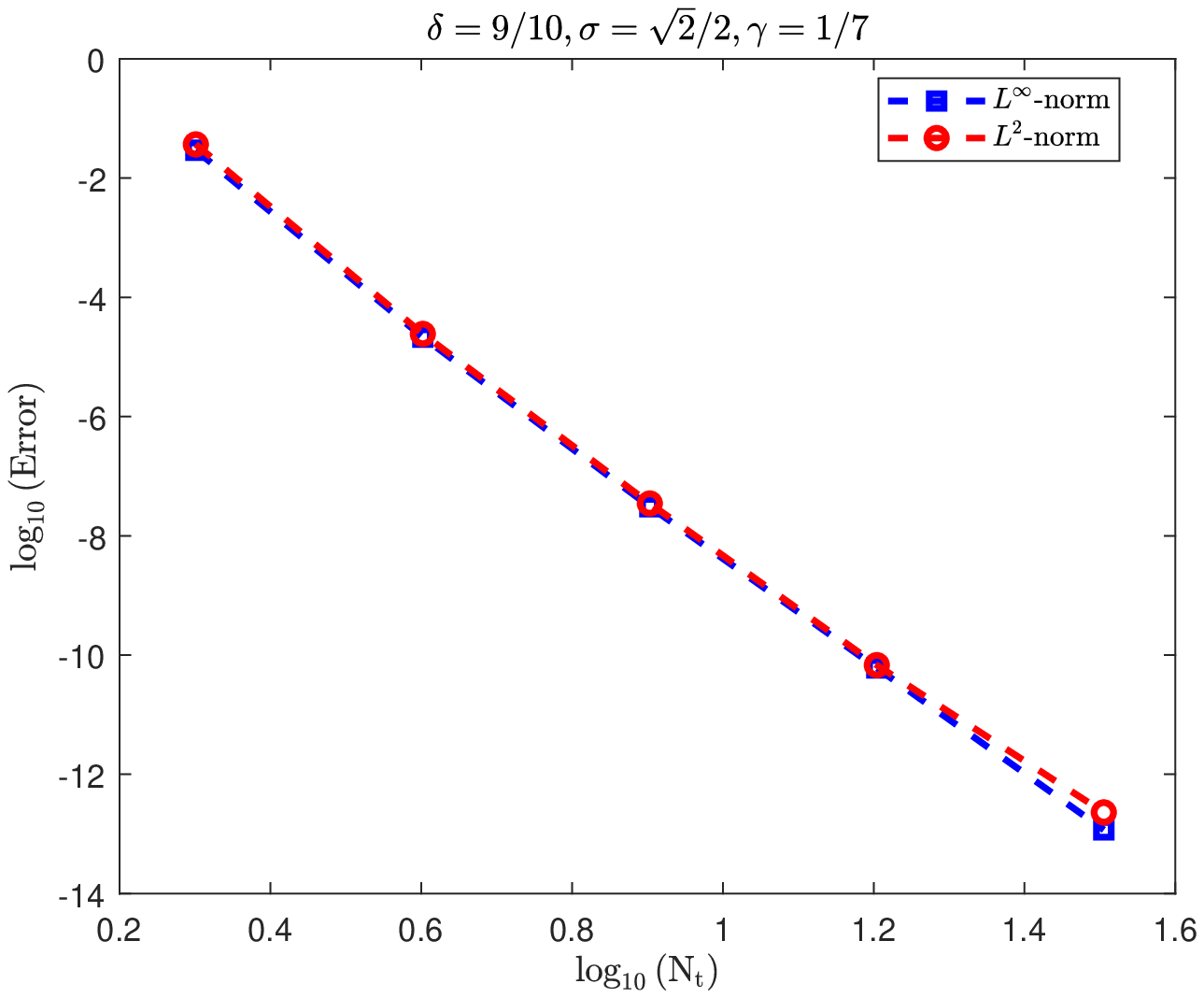}\\
	\caption{ (Example 2) $L^{\infty}$- and $L^2$-errors in log scale versus the polynomial degree for irrational power solution.}\label{ErrPlot2}
\end{figure}

We want to use this example to test the accuracy of the spectral method for nonsmooth solutions.
For the fractional $\sigma=3/5$, we take $\gamma=1/5$ or $1/8.$
The numerical errors versus the polynomial degree $N_t$ for several $\delta$ is plotted in Figure \ref{ErrPlot1}.
It is observed from this figure that the errors decay exponentially as the polynomial degree increases.
For the irrational number $\sigma=\sqrt{2}/2$, we take $\gamma=1/7$.
The obtained result is given in Figure \ref{ErrPlot2}, from which we also observe
the spectral convergence.

\begin{figure}[htbp]
	\centering
	\includegraphics[scale=0.35]{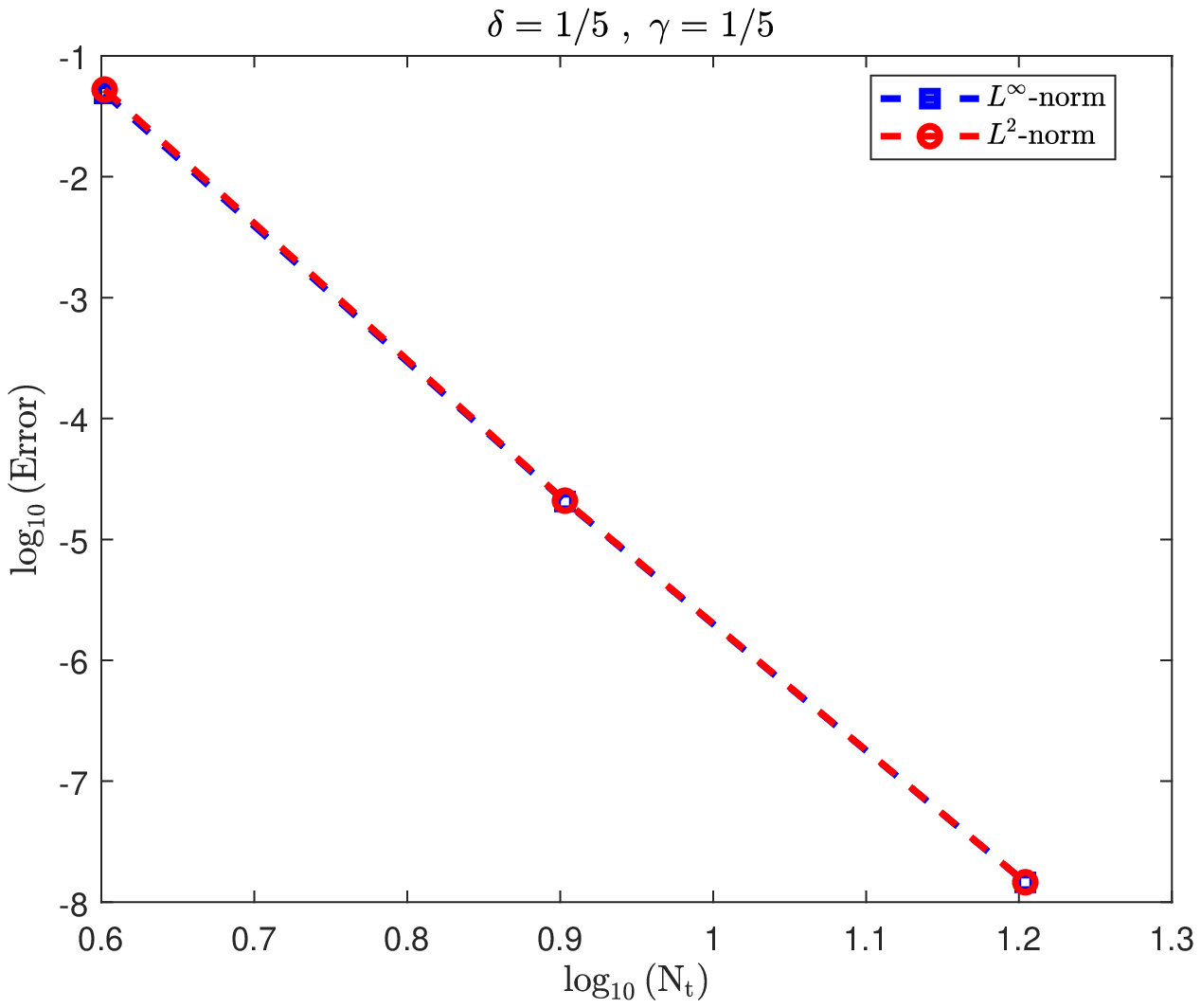}
	\includegraphics[scale=0.35]{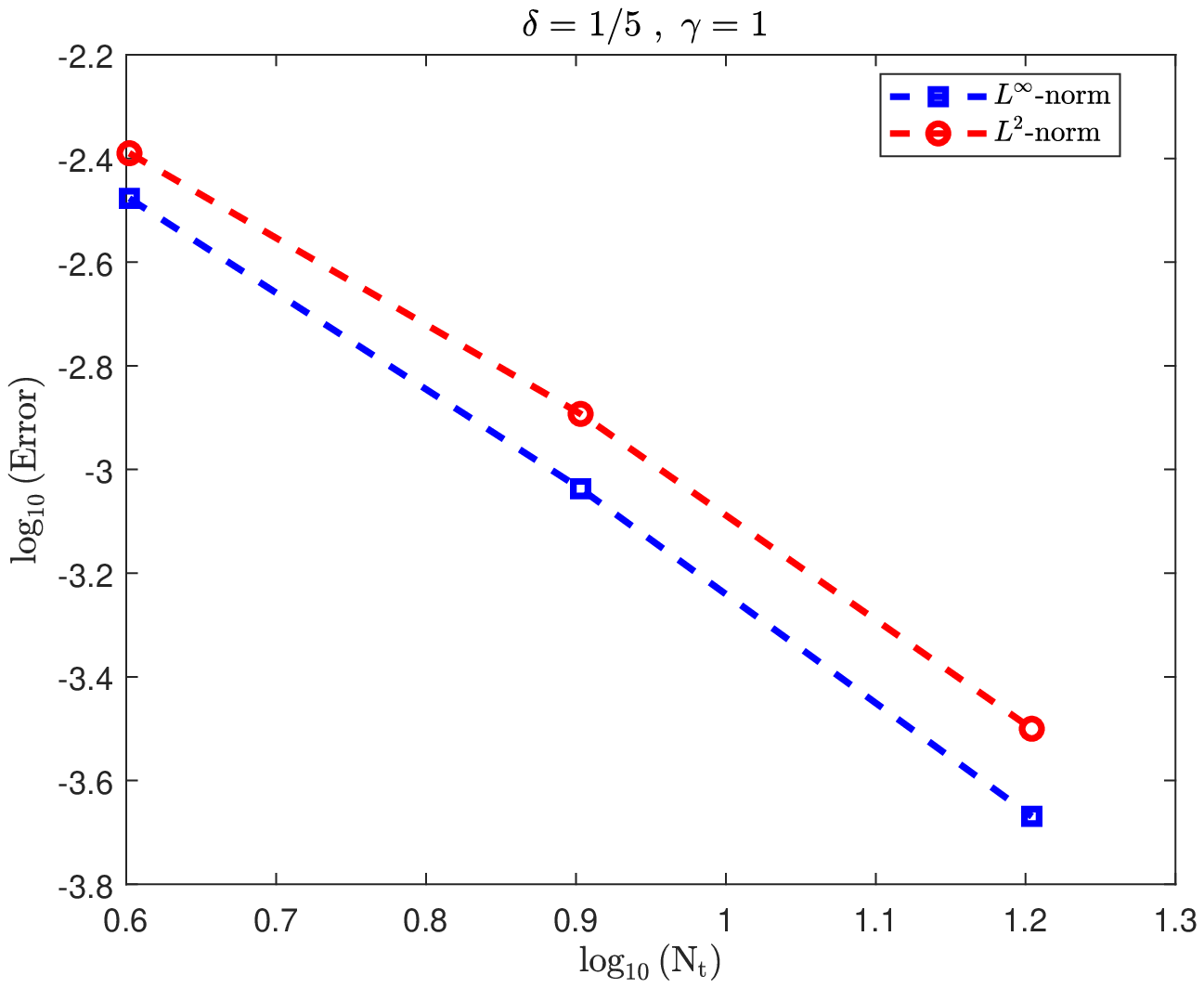}\\
	\includegraphics[scale=0.35]{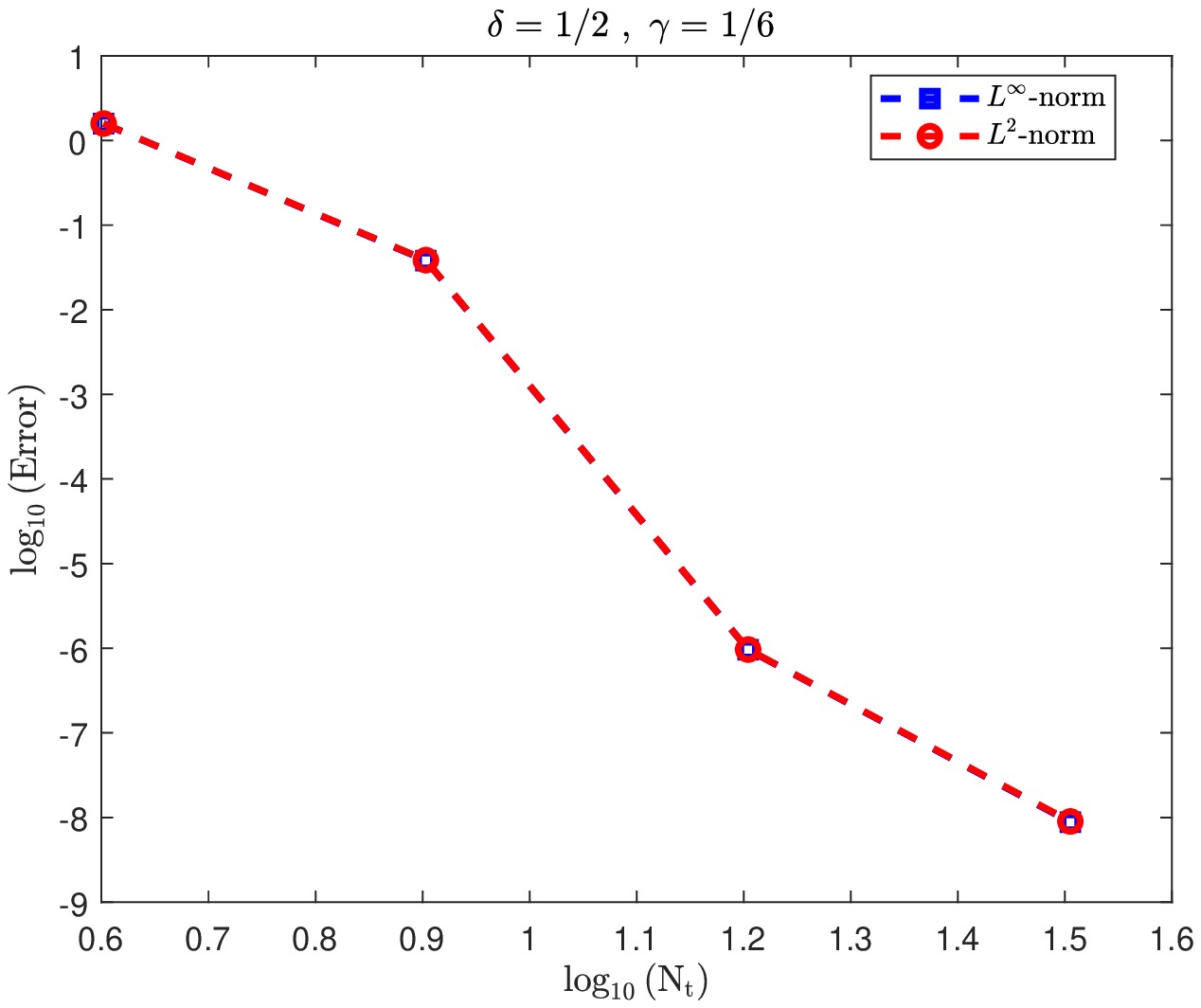}
	\includegraphics[scale=0.35]{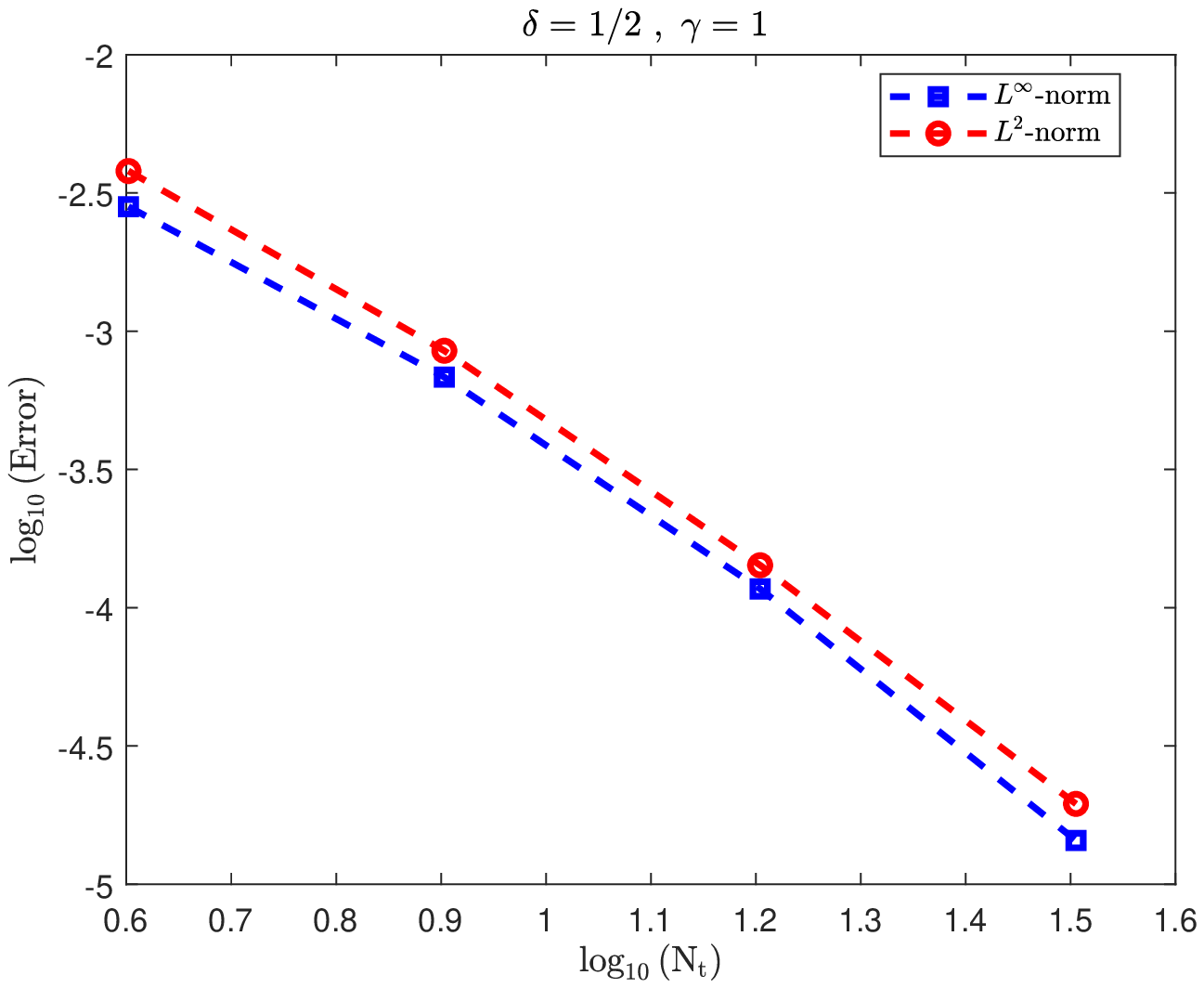}\\
	\includegraphics[scale=0.35]{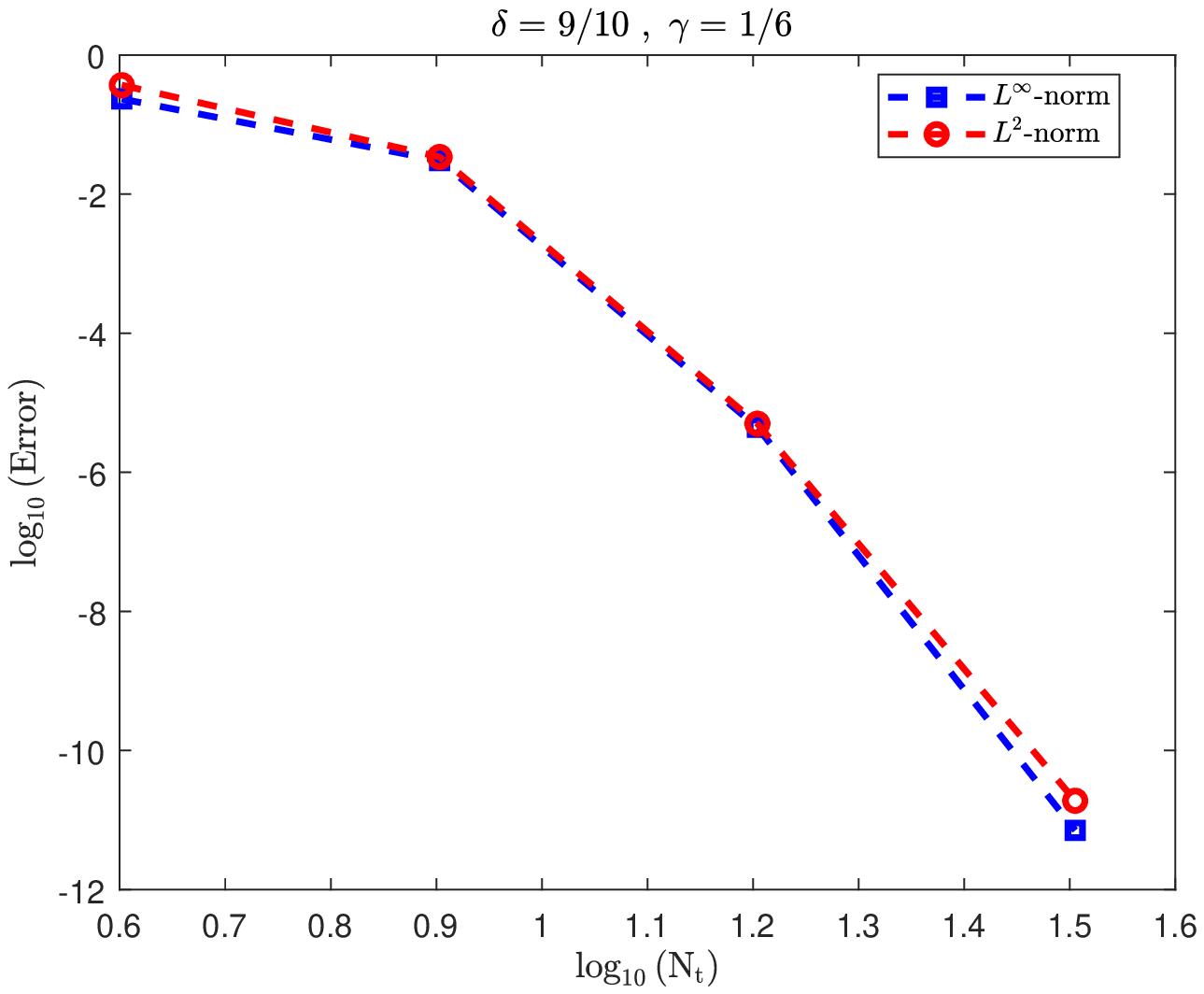}
	\includegraphics[scale=0.35]{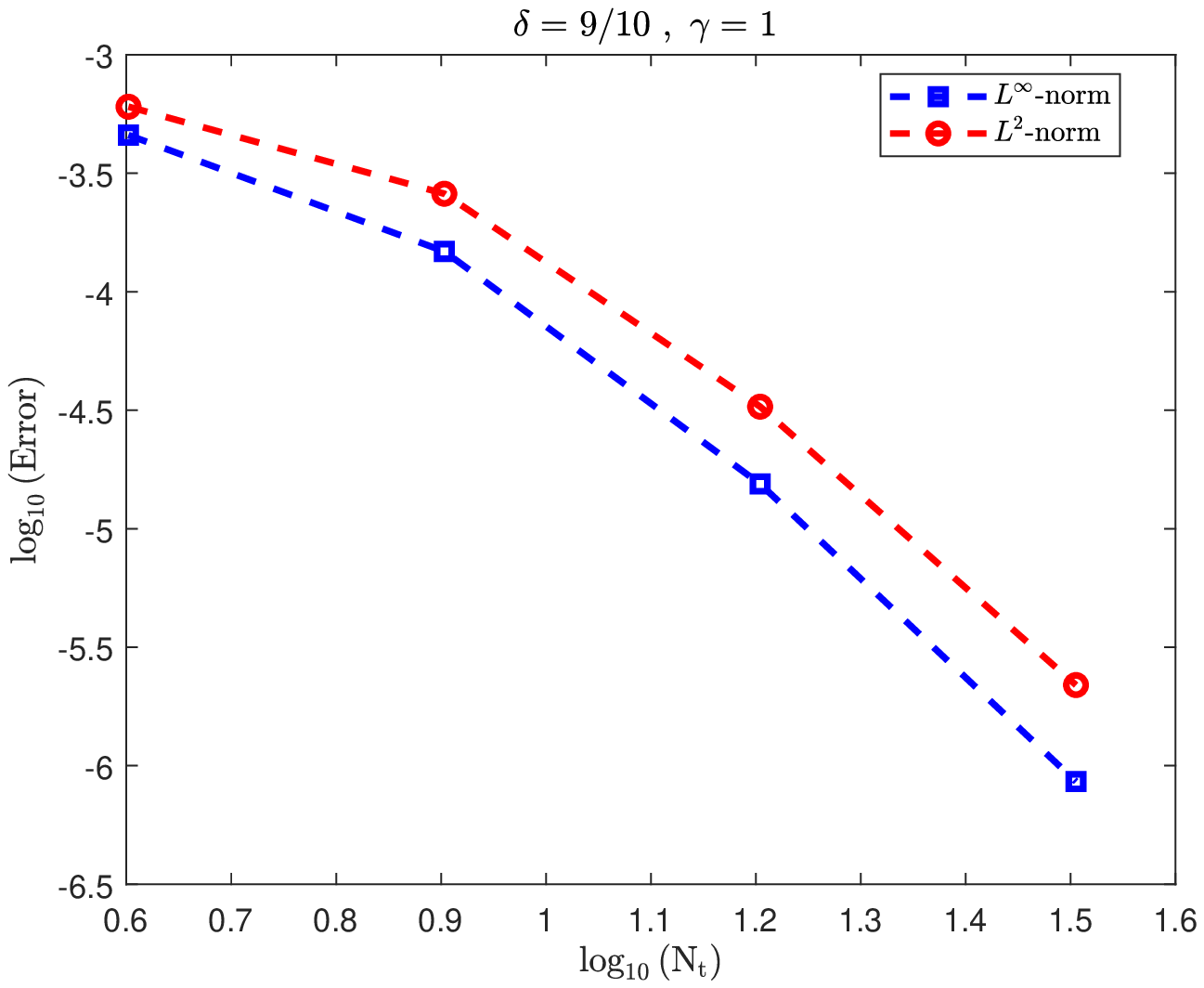}\\
	\caption{(Example 3) $L^{\infty}$- and $L^2$-errors in log scale versus the polynomial degree
$N_t$ for different $\delta$ and $\gamma$.}\label{ErrPlot3}
\end{figure}

\begin{example}(Unknown solution)\label{exam3}
Consider problem \eqref{Caputo_Problem} with a given source function $f(s)=\mathrm{sin}(s)$. In this case the
exact solution and its singularity structure are unknown.
\end{example}
Since the eact solution is unknown, a numerical solution computed with very fine resolution
is served as the reference solution.
The solution qualities are compared for different $\delta$ by two approaches, i.e., our method and usual spectral method, by plotting the
errors versus the polynomial degrees in Figure \ref{ErrPlot3}.
We see that more accurate solutions are obtained by using
$\gamma=1/5$ or $1/6$, compared to the classical spectral method, i.e., $\gamma=1$.

\begin{figure}[htbp]
	\centering
	\includegraphics[scale=0.35]{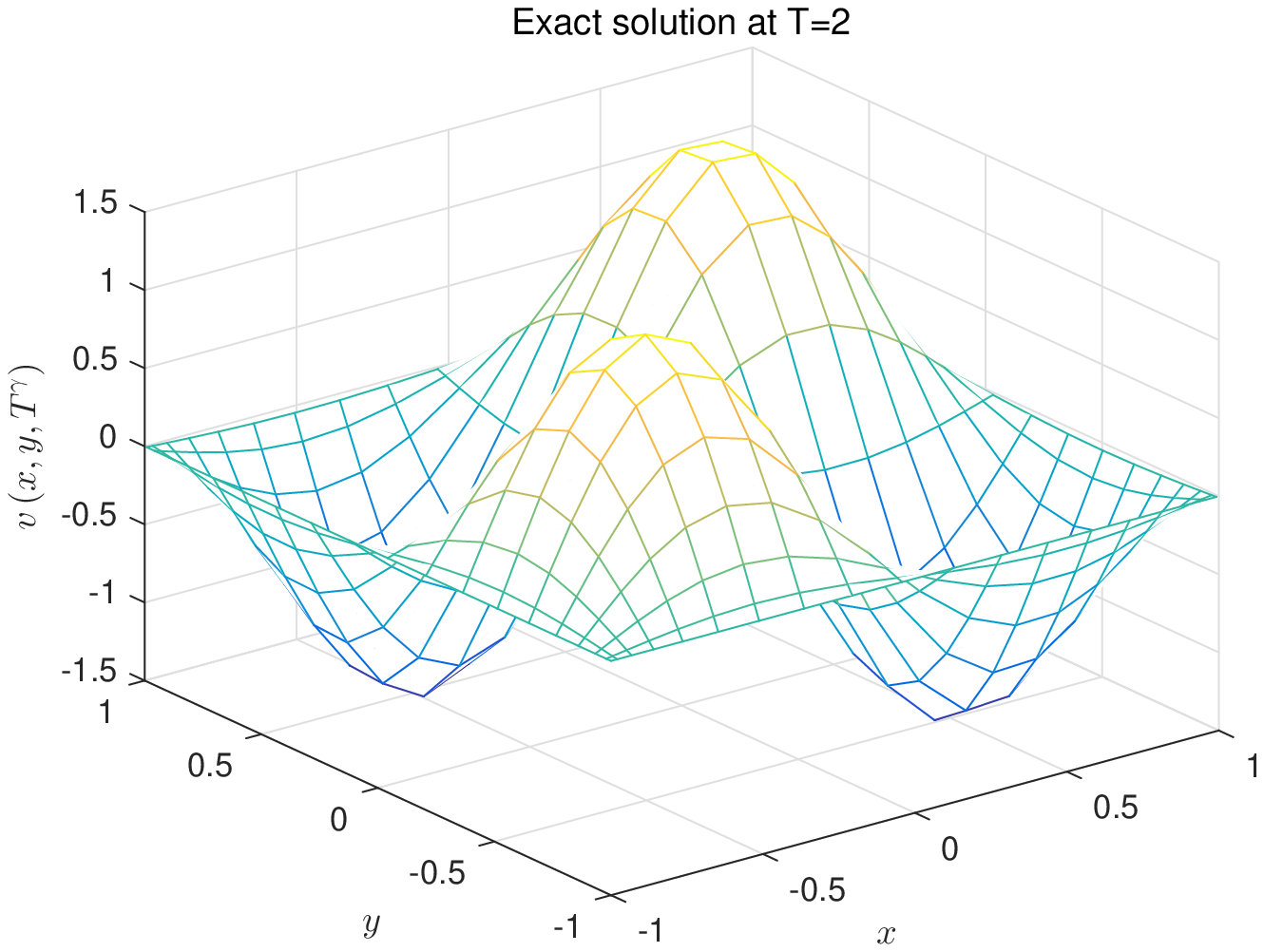}
	\includegraphics[scale=0.35]{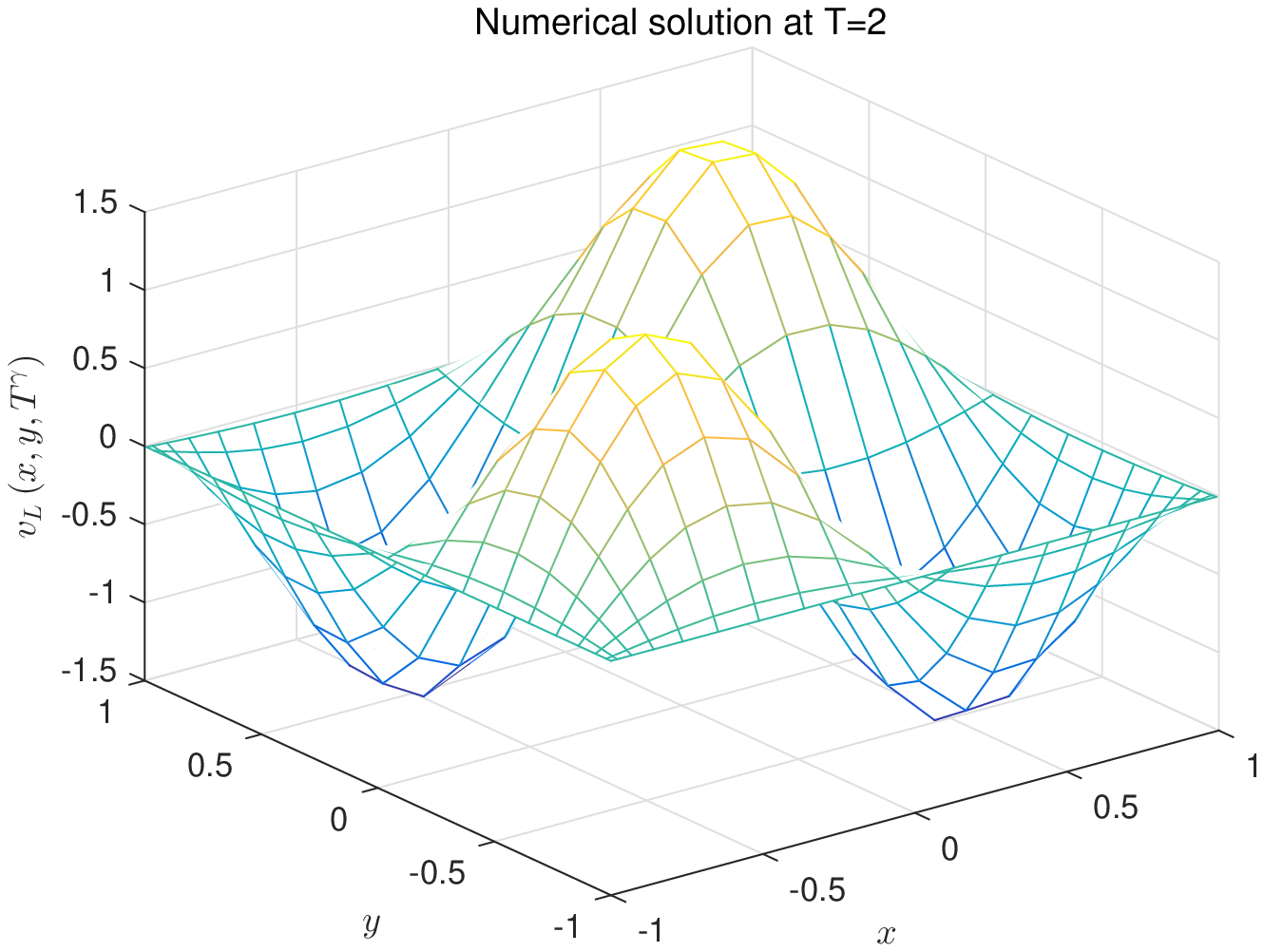}\\
	\includegraphics[scale=0.35]{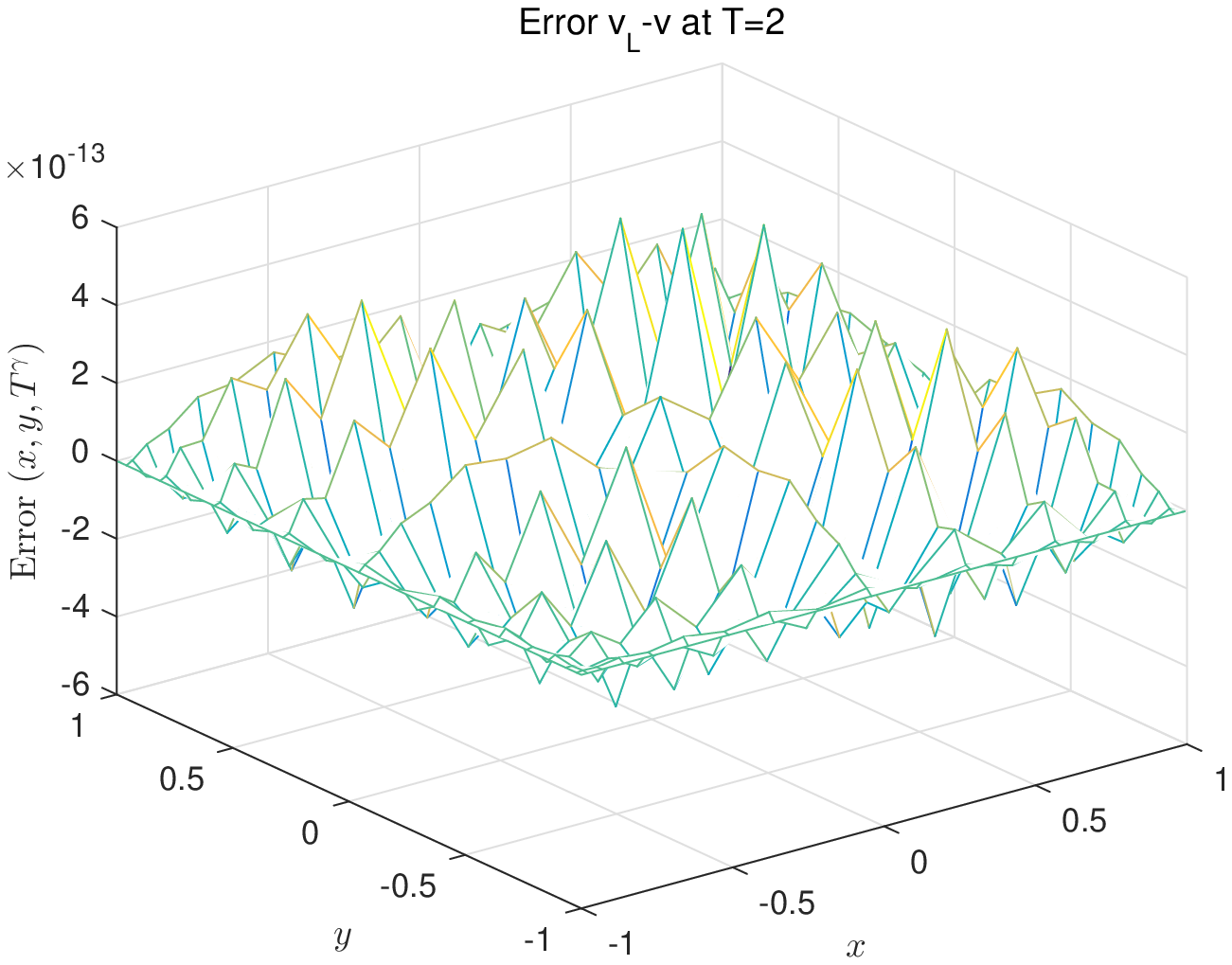}\\
	\caption{(Example 4) The exact solution $v,$ numerical solution $v_L$, and error $v_L-v$ at $T=2$
with $\gamma=1/5$.}\label{ErrPlot7}
\end{figure}

\begin{figure}[htbp]
	\centering
	\includegraphics[scale=0.35]{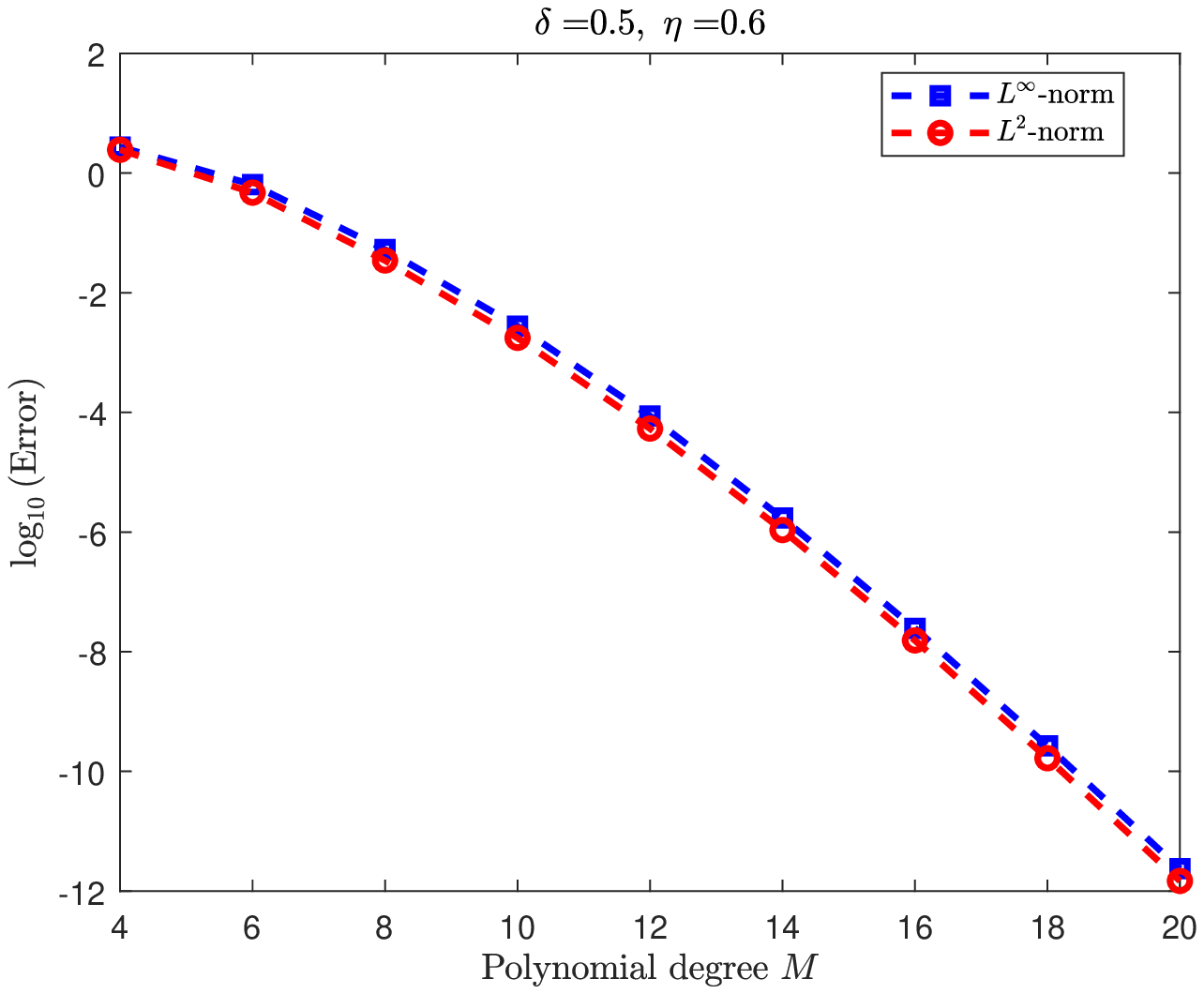}
	\includegraphics[scale=0.35]{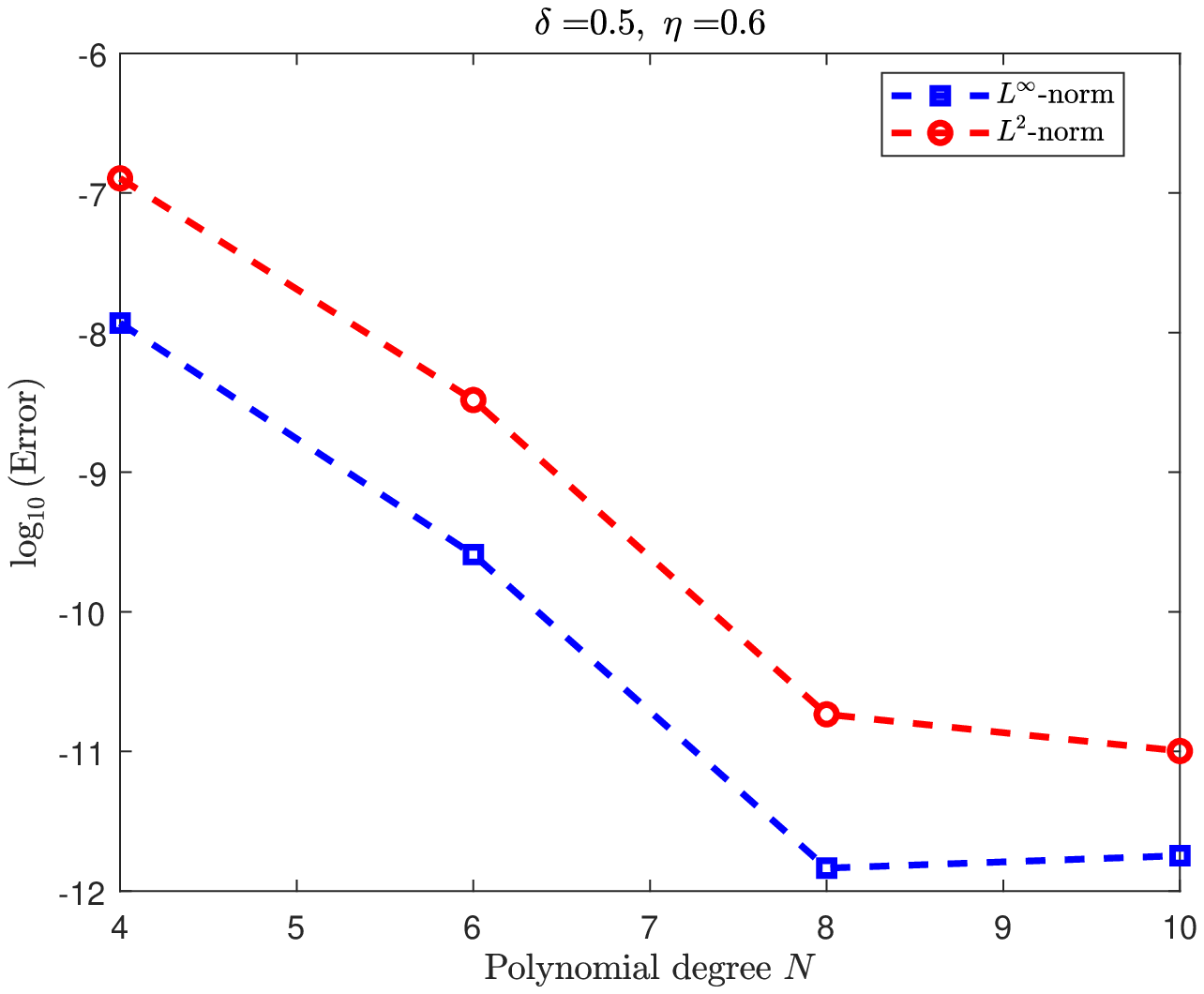}\\
	\caption{(Example 4) $L^{\infty}$- and $L^2$-errors in semi-log scale versus $M$ or $N$ with
	$\delta=0.5,$ $\gamma=1/5$.}\label{ErrPlot8}
\end{figure}

\begin{example}(2D time fractional subdiffusion equation)\label{exam4}
Consider the 2D subdiffusion equation \eqref{TFDE}-\eqref{Bound} with the fabricated exact solution
	$$u(x,y,s)=\mathrm{sin}(\pi x)\mathrm{sin}(\pi y)s^{3/5}.$$
\end{example}
In Figure \ref{ErrPlot7}, we depict the exact solution, numerical solution and error at the final time
computed with the polynomial degree $20$ in both directions.
As shown in this figure, a very accurate solution is obtained with pointwise error
reaching as small as $O(10^{-13})$.
The error history as a function of the polynomial degrees $M$ or $N$, shown
in Figure \ref{ErrPlot8}, confirms the spectral convergence of the used method.


\begin{remark}	
	For the selection of parameter $\gamma,$ our fundamental principle is to make $v(\cdot)=u((\cdot)^{1/\gamma})$ sufficiently smooth which can be made according to the following strategy:
	
	Case I: if the solution $u$ is smooth, the optimal value is $\gamma=1;$
	
	Case II: if the source term $g(\cdot)$ is smooth, then
	$(a)$ when $\delta$ is a rational number $p/q$, the best choice is $\gamma=1/q$. Theoretically $\gamma=1/nq~(n=1,2,\cdots)$ works too, but larger $n$ leads to larger amount of calculation;
	$(b)$ when $\delta$ is an irrational number, there is no suitable value of $\gamma$ to make $u((\cdot)^{1/\gamma})$ smooth. In this case, we can take $\gamma=1/q$ with a reasonably large $q$ such that $u((\cdot)^{1/\gamma})$ is smooth enough.
\end{remark}

\section{Concluding remarks}\label{Sect6}

A novel spectral method has been proposed and analyzed for the subdiffusion equation.
The main novelty of the proposed method is its variational framework based on
fractional Sobolev spaces.
The idea was to first apply suitable variable transformation
to re-scale the underlying equation, then construct spectral methods for the re-scaled
equation. This is particularly useful in numerical solutions of fractional differential equations,
to which the solution is often singular and can be smoothed by using appropriate transformation.
For this purpose, a new variational framework was established based on the fractional Sobolev spaces,
which allows constructing and analyzing numerical methods
following the standard Galerkin approach.
Our theoretical and numerical investigation showed that the proposed method using suitable  transformation is exponentially convergent for general right hand side functions, even though the exact solution has limited regularity.
Implementation details was also provided, along with
a series of numerical examples to demonstrate the efficiency of the proposed method.

It is worthy to mention here a number of points: First, with some specific choices of the transformation
function, the new method can be proved to be equivalent to the M\"untz spectral method, recently proposed in a series of
papers \cite{Hou2017A, Hou2018Muntz}. The latter was based on the
M\"untz polynomial approximation to the original equation;
Secondly, although the error analysis was carried out only for a particular transformation, it seems
extendable to some other choices;
Finally, compared to the M\"untz spectral method, the main benefit of the current method
may be its flexibility in choosing the transformation function. This makes the new method
applicable to a larger class of problems.

\section*{Acknowledgements}
This research is supported by the NSFC grant 11971408.
%

\bibliographystyle{elsarticle-num}
\bibliography{refs}
\end{document}